\documentclass{amsproc}

\usepackage{amsmath}
\usepackage{amssymb}
\usepackage{amsxtra}
\usepackage{latexsym}
\usepackage{ifthen}

\usepackage[dvips]{graphicx}
\usepackage{epstopdf}
\usepackage{epsfig}
\usepackage{epic}
\usepackage{eepic}

\def\Xint#1{\mathchoice
   {\XXint\displaystyle\textstyle{#1}}
   {\XXint\textstyle\scriptstyle{#1}}
   {\XXint\scriptstyle\scriptscriptstyle{#1}}
   {\XXint\scriptscriptstyle\scriptscriptstyle{#1}}
   \!\int}
\def\XXint#1#2#3{{\setbox0=\hbox{$#1{#2#3}{\int}$}
     \vcenter{\hbox{$#2#3$}}\kern-.5\wd0}}
\def\dashint{\Xint-}

\newtheorem{theorem}{Theorem}[section]
\newtheorem{lemma}[theorem]{Lemma}
\newtheorem{corollary}[theorem]{Corollary}
\newtheorem{propo}[theorem]{Proposition}

\theoremstyle{definition}

\theoremstyle{remark}
\newtheorem{remark}[theorem]{Remark}

\numberwithin{equation}{section}

\def\Xint#1{\mathchoice
   {\XXint\displaystyle\textstyle{#1}}%
   {\XXint\textstyle\scriptstyle{#1}}%
   {\XXint\scriptstyle\scriptscriptstyle{#1}}%
   {\XXint\scriptscriptstyle\scriptscriptstyle{#1}}%
   \!\int}
\def\XXint#1#2#3{{\setbox0=\hbox{$#1{#2#3}{\int}$}
     \vcenter{\hbox{$#2#3$}}\kern-.5\wd0}}
\def\dashint{\Xint-}

\begin{document}

\title[Ring homeomorphisms and prime ends]{RING HOMEOMORPHISMS AND PRIME ENDS}

\author{Vladimir Gutlyanskii}
\address{Institute of Applied Mathematics and Mechanics, National Academy of Sciences of Ukraine,
 Slavyansk, Ukraine}
\email{vladimirgut@mail.ru, vladimirgutlyanskii@yahoo.com}

\author{Vladimir Ryazanov}
\address{Institute of Applied Mathematics and Mechanics, National Academy of Sciences of Ukraine,
 Slavyansk, Ukraine}
\email{vlryazanov1@rambler.ru, vl$\underline{\ \
}$\,ryazanov1@mail.ru}

\author{Eduard Yakubov}
\address{Holon Institute of Technology, Holon, Israel}
\email{yakubov@hit.ac.il , eduardyakubov@gmail.com}

\date{\today}



\begin{abstract}

We show that every homeomorphic $W^{1,1}_{\rm loc}$ solution $f$ of
a Beltrami equation $\overline{\partial}f=\mu\,\partial f$ in a
domain $D\subseteq\Bbb C$ is the so--called ring
$Q-$ho\-meo\-mor\-phism with $Q(z)=K^T_{\mu}(z, z_0)$ where
$K^T_{\mu}(z, z_0)$ is the tangent (angular) dilatation quotient of
the equation with respect to an arbitrary point $z_0\in
{\overline{D}}$. In this connection, we develop the theory of the
boundary behavior of the ring $Q-$homeo\-mor\-phisms with respect to
prime ends. On this basis, we show that, for wide classes of
degenerate Beltrami equations $\overline{\partial}f=\mu\,\partial
f$, there exist regular solutions of the Dirichlet problem in
arbitrary simply connected domains in $\Bbb C$ and pseudoregular and
multivalent solutions in arbitrary finitely  connected domains in
$\Bbb C$ with boundary datum $\varphi$ that are continuous with
respect to the topology of prime ends.
\end{abstract}

\maketitle

\tableofcontents

\newpage

\section{Introduction}

Let $D$ be a domain in the complex plane ${\Bbb C}$ and
$\mu:D\to{\Bbb C}$ be a measurable function with $|\mu(z)|<1$ a.e.
(almost everywhere) in $D$. A {\bf Beltrami equation} is an equation
of the form
\begin{equation}\label{eqBeltrami} f_{\bar z}=\mu(z)\,f_z\end{equation} where
$f_{\bar z}=\overline{\partial}f=(f_x+if_y)/2$,
$f_{z}=\partial f=(f_x-if_y)/2$, $z=x+iy$, and $f_x$ and $f_y$ are
partial derivatives of $f$ in $x$ and $y$, correspondingly.

\medskip

Boundary value problems for the Beltrami equations are due to the
well-known Riemann dissertation in the case of $\mu(z)=0$ and to the
papers of Hilbert (1904, 1924) and Poincare (1910) for the
corresponding Cauchy--Riemann system.

\medskip

The classic {\bf Dirichlet problem} for the Beltrami equation
(\ref{eqBeltrami}) in a domain $D\subset{\Bbb C}$ is the problem on
the existence of a continuous function $f:D\to{\Bbb C}$ having
partial derivatives of the first order a.e., satisfying
(\ref{eqBeltrami}) a.e. and such that
\begin{equation}\label{eqGrUsl}\lim\limits_{z\to\zeta}{\rm
Re}\,f(z)=\varphi(\zeta)\qquad\forall\ \zeta\in\partial
D\end{equation} for a prescribed continuous function
$\varphi:\partial D\to{\Bbb R},$ see, e.g., \cite{Boj} and
\cite{Vekua}.

\medskip

The function $\mu$ is called the {\bf complex coefficient} and
\begin{equation}\label{eqKPRS1.1}K_{\mu}(z)=\frac{1+|\mu(z)|}{1-|\mu(z)|}\end{equation}
the {\bf dilatation quotient} of the equation (\ref{eqBeltrami}).
The Beltrami equation (\ref{eqBeltrami}) is said to be {\bf
degenerate} if ${\rm ess}\,{\rm sup}\,K_{\mu}(z)=\infty$.

\medskip

The existence of homeomorphic $W^{1,1}_{\rm loc}$ solutions was
recently established to many degenerate Beltrami equations, see,
e.g., \cite{GRSY} and \cite{MRSY} and references therein. The theory
of boundary behavior of homeomorphic solutions and of the Dirichlet
problem for a wide circle of degenerate Beltrami equations in Jordan
domains was developed in \cite{KPR}, \cite{KPRS}, \cite{KPRS_*} and
\cite{RSSY}.

\medskip

To study the similar problem in domains with more complicated
boundaries we need to apply the theory of prime ends by
Caratheodory, see, e.g., his paper \cite{Car$_2$} or Chapter 9 in
monograph \cite{CL}.

\medskip

The main difference in the case is that $\varphi$ should be already
a function of a boundary element (prime end $P$) but not of a
boundary point. Moreover, (\ref{eqGrUsl}) should be replaced by the
condition
\begin{equation}\label{eq1002P} \lim\limits_{n\to\infty}\ {\rm
Re}\, f(z_n)\ =\ \varphi(P)
\end{equation} for all sequences of points $z_n\in D$ converging to
prime ends  $P$ of the domain $D$. Note that (\ref{eq1002P}) is
equivalent to the condition that
\begin{equation}\label{eq1002PP} \lim\limits_{z\to P}\ {\rm
Re}\, f(z)\ =\ \varphi(P)
\end{equation} along any ways in $D$ going to the prime ends $P$ of
the domain $D$.

\medskip

Later on, $\overline{D}_P$ denotes the completion of the domain $D$
by its prime ends and $E_D$ denotes the space of these prime ends,
both with the topology of prime ends described in Section 9.5 of
monograph \cite{CL}. In addition, continuity of mappings
$f:\overline{D}_P\to\overline{D^{\prime}}_P$ and boundary functions
$\varphi : E_D\to\Bbb R$ should mean with respect to the given
topology of prime ends.

\begin{remark}\label{METRIC} This topology can be described in terms
of metrics. Namely, as known, every bounded finitely connected
domain $D$ in ${\Bbb C}$ can be mapped by a conformal mapping $g_0$
onto the so-called circular domain $D_0$ whose boun\-da\-ry consists
of a finite collection of mutually disjoint circles and isolated
points, see, e.g., Theorem V.6.2 in \cite{Goluzin}. Moreover,
isolated singular points of bounded conformal mappings are removable
by Theorem 1.2 in \cite{CL} due to Weierstrass. Hence isolated
points of $\partial D$ correspond to isolated points of $\partial
D_0$ and inversely.


Reducing this case to the Caratheodory theorem, see, e.g., Theorem
9.4 in \cite{CL} for simple connected domains, we have a natural
one-to-one correspondence between points of $\partial D_0$ and prime
ends of the domain $D$. Determine in $\overline{D}_P$ the metric
$\rho_0(p_1,p_2)=\left|{\widetilde
{g_0}}(p_1)-{\widetilde{g_0}}(p_2)\right|$ where ${\widetilde
{g_0}}$ is the extension of $g_0$ to $\overline{D}_P$ just
mentioned.


If $g_*$ is another conformal mapping of the domain $D$ on a
circular domain $D_*$, then the corresponding metric
$\rho_*(p_1,p_2)=\left|{\widetilde{g_*}}(p_1)-{\widetilde{g_*}}(p_2)\right|$
generates the same convergence in $\overline{D}_P$ as the metric
$\rho_0$ because $g_0\circ g_*^{-1}$ is a conformal mapping between
the domains $D_*$ and $D_0$ that is extended to a homeomorphism
between $\overline{D_*}$ and $\overline{D_0}$, see, e.g., Theorem
V.6.1$^{\prime}$ in \cite{Goluzin}. Con\-se\-quent\-ly, the given
metrics induce the same topology in the space $\overline{D}_P$.


This topology coincides with topology of prime ends described in
inner terms of the domain $D$ in Section 9.5 of \cite{CL}. Later on,
we prefer to apply the description of the topology of prime ends in
terms of the given metrics because it is more clear, more convenient
and it is important for us just metrizability of $\overline{D}_P$.
Note also that the space $\overline{D}_P$ for every bounded finitely
connected domain $D$ in ${\Bbb C}$ with the given topology is
compact because the closure of the circular domain $D_0$ is a
compact space and by the construction $\widetilde
{g_0}:\overline{D}_P\to{\overline {D_0}}$ is a homeomorphism.
\end{remark}

\medskip

Applying the description of the topology of prime ends in Section
9.5 of \cite{CL}, we reduce the case of bounded finitely connected
domains to Theorem 9.3 in \cite{CL} for simple connected domains and
obtain the following useful fact.

\medskip

\begin{lemma}\label{thabc3} {\it Each prime end $P$ of a bounded finitely
connected domain $D$ in $\Bbb C$ contains a chain of cross--cuts
$\sigma_m$ lying on circles $S(z_0,r_m)$ with $z_0\in\partial D$ and
$r_m\to0$ as $m\to\infty$.}
\end{lemma}

\medskip

Given a point $z_0$ in $\Bbb C$, we apply here the more refined
quantity than $K_{\mu}(z)\ $:
\begin{equation}\label{eqTangent} K^T_{\mu}(z,z_0)\ =\
\frac{\left|1-\frac{\overline{z-z_0}}{z-z_0}\ \mu
(z)\right|^2}{1-|\mu (z)|^2} \end{equation} that is called the {\bf
tangent dilatation quotient} of the Beltrami equation
(\ref{eqBeltrami}) with respect to $z_0$, see, e.g., \cite{RSY$_3$},
cf. the corresponding terms and notations in \cite{And3, GMSV$_2$,
Le} and \cite{RW}. The given term was first introduced in
\cite{RSY$_3$} and its geometric sense was described in \cite{RSSY},
see also \cite{MRSY}, Section 11.3. Note that
\begin{equation}\label{eqConnect} K^{-1}_{\mu}(z)\ \leqslant\ K^T_{\mu}(z,z_0)\
\leqslant\ K_{\mu}(z)\ \ \ \ \ \ \ \ \ \ \forall\ z\in D\ \ \
\forall\ z_0\in \Bbb C
\end{equation} and the given estimates are precise. The quantity (\ref{eqTangent})
takes into account not only the modulus of the complex coefficient
$\mu$ but also its argument.

\medskip

Throughout this paper, $B(z_0,r)=\{z\in{\Bbb C}:|z-z_0|<r\}$, ${\Bbb
D}=B(0,1)$, $S(z_0,r)=\{z\in{\Bbb C}:|z-z_0|=r\}$,
$R(z_0,r_1,r_2)=\{z\in{\Bbb C}:r_1<|z-z_0|<r_2\}$.

\bigskip

\section{Regular domains}

First of all, recall the following topological notion. A domain
$D\subset{\Bbb C}$ is said to be {\bf locally connected at a point}
$z_0\in\partial D$ if, for every neighborhood $U$ of the point
$z_0$, there is a neighborhood $V\subseteq U$ of $z_0$ such that
$V\cap D$ is connected. If this condition holds for all $z_0\in
\partial D$, then $D$ is said to be locally connected on $\partial
D$. For a domain that is locally connected on $\partial D$, there is
a natural one--to--one correspondence between prime ends of $D$ and
points of $\partial D$ and the topology of prime ends coincides with
the Euclidean topology. Note that every Jordan domain $D$ in ${\Bbb
C}$ is locally connected on  $\partial D$, see, e.g., \cite{Wi}, p.
66.

\medskip

Now, recall that the {\bf (conformal) modulus} of a family $\Gamma$
of paths $\gamma$ in ${\Bbb C}$ is the quantity
\begin{equation}\label{eqModul} M(\Gamma)=\inf_{\varrho\in{\rm adm}\,\Gamma}\int\limits_{\Bbb
C}\varrho^2(z)\,dm(z)\end{equation} where a Borel function
$\varrho:{\Bbb C}\to[0,\infty]$ is {\bf admissible} for $\Gamma$,
write $\varrho\in{\rm adm}\,\Gamma$, if
\begin{equation}\label{eqAdm}\int\limits_{\gamma}\varrho\,ds\geqslant1\quad\forall\ \gamma\in\Gamma\,.\end{equation}
Here $s$ is a natural parameter of the are length on $\gamma$.

\medskip

Later on, given sets $A,$ $B$ and $C$ in $\mathbb{C},$ $\Delta (A,
B; C)$ denotes a family of all paths $\gamma : [a, b] \to
\mathbb{C}$ joining $A $ and $B$ in $C,$ i.e. $\gamma (a) \in A,$
$\gamma (b) \in B$ and $\gamma (t) \in C$ for all $t \in (a,b).$

\medskip

We say that $\partial D$ is {\bf weakly flat at a point}
$z_0\in\partial D$ if, for every neighborhood $U$ of the point
$z_0$ and every number $P>0$, there is a neighborhood $V\subset U$
of $z_0$ such that \begin{equation}\label{eq1.5KR}
M(\Delta(E,F;D))\geqslant P\end{equation} for all continua $E$ and
$F$ in $D$ intersecting $\partial U$ and $\partial V$. We say that
$\partial D$ is {\bf weakly flat} if it is weakly flat at each
point $z_0\in\partial D$.

\medskip

We also say that a point $z_0\in\partial D$ is {\bf strongly
accessible} if, for every neigh\-bor\-hood $U$ of the point $z_0$,
there exist a compactum $E$ in $D$, a neigh\-bor\-hood $V\subset U$
of $z_0$ and a number $\delta>0$ such that
\begin{equation}\label{eq1.6KR}M(\Delta(E,F;D))\geqslant\delta\end{equation}
for all continua $F$ in $D$ intersecting $\partial U$ and
$\partial V$. We say that $\partial D$ is {\bf strongly
accessible} if each point $z_0\in\partial D$ is strongly
accessible.

\medskip

It is easy to see that if a domain $D$ in ${\Bbb C}$ is weakly flat
at a point $z_0\in\partial D$, then the point $z_0$ is strongly
accessible from $D$. The following fact is fundamental, see, e.g.,
Lemma 5.1 in \cite{KR$_2$} or Lemma 3.15 in \cite{MRSY}.


\begin{lemma}\label{lem3.68} If a domain $D$ in $\Bbb C$ is weakly flat
at a point $z_0\in\partial D$, then $D$ is locally connected at
$z_0$. \end{lemma}


The notions of strong accessibility and weak flatness at boundary
points of a domain in ${\Bbb C}$ defined in \cite{KR$_0$}, see also
\cite{KR$_2$} and \cite{RS}, are localizations and
ge\-ne\-ra\-li\-za\-tions of the corresponding notions introduced in
\cite{MRSY$_5$} and \cite{MRSY$_6$}, cf. with the properties $P_1$
and $P_2$ by V\"ais\"al\"a in \cite{Va} and also with the
quasiconformal accessibility and the quasiconformal flatness by
N\"akki in \cite{Na$_1$}.

\medskip

A domain $D\subset{\Bbb C}$ is called a {\bf quasiextremal distance
domain}, abbr. {\bf QED-domain}, see \cite{GM}, if
\begin{equation}\label{e:7.1}M(\Delta(E,F;{\Bbb C})\leqslant K\cdot
M(\Delta(E,F;D))\end{equation} for some $K\geqslant1$ and all pairs
of nonintersecting continua $E$ and $F$ in $D$.

\medskip

It is well known, see, e.g., Theorem 10.12 in \cite{Va}, that
\begin{equation}\label{eqKPR2.2}M(\Delta(E,F;{\Bbb C}))\geqslant\frac{2}{\pi}\log{\frac{R}{r}}\end{equation}
for any sets $E$ and $F$ in ${\Bbb C}$ intersecting all the
circles $S(z_0,\rho)$, $\rho\in(r,R)$. Hence a QED-domain has a
weakly flat boundary. One example in \cite{MRSY}, Section 3.8,
shows that the inverse conclusion is not true even in the case of
simply connected domains in ${\Bbb C}$.

\medskip

A domain $D\subset{\Bbb C}$ is called a {\bf uniform domain} if
each pair of points $z_1$ and $z_2\in D$ can be joined with a
rectifiable curve $\gamma$ in $D$ such that
\begin{equation}\label{e:7.2} s(\gamma)\leqslant
a\cdot|z_1-z_2|\end{equation} and \begin{equation}\label{e:7.3}
\min\limits_{i=1,2}\ s(\gamma(z_i,z))\leqslant b\cdot{\rm
dist}(z,\partial D)\end{equation} for all $z\in\gamma$ where
$\gamma(z_i,z)$ is the portion of $\gamma$ bounded by $z_i$ and
$z$, see \cite{MaSa}. It is known that every uniform domain is a
QED-domain but there exist QED-domains that are not uniform, see
\cite{GM}. Bounded convex domains and bounded domains with smooth
boundaries are simple examples of uniform domains and,
consequently, QED-domains as well as domains with weakly flat
boundaries.

\medskip

It is also often met with the so-called Lipschitz domains in the
mapping theory and in the theory of differential equations. Recall
first that $\varphi:U\to{\Bbb C}$ is said to be a {\bf Lipcshitz
map} provided $|\varphi(z_1)-\varphi(z_2)|\leqslant M\cdot|z_1-z_2|$
for some $M<\infty$ and for all $z_1$ and $z_2\in U$, and a {\bf
bi-Lipcshitz map} if in addition
$M^*|z_1-z_2|\leqslant|\varphi(z_1)-\varphi(z_2)|$ for some $M^*>0$
and for all $z_1$ and $z_2\in U$. They say that $D$ in $\Bbb{C}$ is
a {\bf Lipschitz domain} if every point $z_0\in\partial D$ has a
neighborhood $U$ that can be mapped by a bi-Lipschitz homeomorphism
$\varphi$ onto the unit disk ${\Bbb D}$ in ${\Bbb C}$ such that
$\varphi(\partial D\cap U)$ is the intersection of ${\Bbb D}$ with
the real axis. Note that a bi-Lipschitz homeomorphism is
quasiconformal and, consequently, the modulus is quasiinvariant
under such a mapping. Hence the Lipschitz domains have weakly flat
boundaries.

\section{BMO, VMO and FMO functions}

Recall that a real-valued function $u$ in a domain $D$ in ${\Bbb C}$
is said to be of {\bf bounded mean oscillation} in $D$, abbr.
$u\in{\rm BMO}(D)$, if $u\in L_{\rm loc}^1(D)$ and
\begin{equation}\label{lasibm_2.2_1}\Vert u\Vert_{*}:=
\sup\limits_{B}{\frac{1}{|B|}}\int\limits_{B}|u(z)-u_{B}|\,dm(z)<\infty\,,\end{equation}
where the supremum is taken over all discs $B$ in $D$, $dm(z)$
corresponds to the Lebesgue measure in ${\Bbb C}$ and
$$u_{B}={\frac{1}{|B|}}\int\limits_{B}u(z)\,dm(z)\,.$$ We write $u\in{\rm BMO}_{\rm loc}(D)$ if
$u\in{\rm BMO}(U)$ for every relatively compact subdomain $U$ of $D$
(we also write BMO or ${\rm BMO}_{\rm loc }$ if it is clear from the
context what $D$ is).

\medskip

The class BMO was introduced by John and Nirenberg (1961) in the
paper \cite{JN} and soon became an important concept in harmonic
analysis, partial differential equations and related areas; see,
e.g., \cite{HKM} and \cite{RR}.
\medskip

A function $\varphi$ in BMO is said to have {\bf vanishing mean
oscillation}, abbr. $\varphi\in{\rm VMO}$, if the supremum in
(\ref{lasibm_2.2_1}) taken over all balls $B$ in $D$ with
$|B|<\varepsilon$ converges to $0$ as $\varepsilon\to0$. VMO has
been introduced by Sarason in \cite{Sarason}. There are a number of
papers devoted to the study of partial differential equations with
coefficients of the class VMO, see, e.g., \cite{CFL, ISbord,
MRV$^*$, Pal} and \cite{Ra}.

\medskip

\begin{remark}\label{rem1} Note that
$W^{\,1,2}\left({{D}}\right) \subset VMO \left({{D}}\right),$ see,
e.g., \cite{BN}.
\end{remark}

\medskip

Following \cite{IR}, we say that a function $\varphi:D\to{\Bbb R}$
has {\bf finite mean oscillation} at a point $z_0\in D$, abbr.
$\varphi\in{\rm FMO}(z_0)$, if
\begin{equation}\label{FMO_eq2.4}\overline{\lim\limits_{\varepsilon\to0}}\ \ \
\dashint_{B(z_0,\varepsilon)}|{\varphi}(z)-\widetilde{\varphi}_{\varepsilon}(z_0)|\,dm(z)<\infty\,,\end{equation}
where \begin{equation}\label{FMO_eq2.5}
\widetilde{\varphi}_{\varepsilon}(z_0)=\dashint_{B(z_0,\varepsilon)}
{\varphi}(z)\,dm(z)\end{equation} is the mean value of the function
${\varphi}(z)$ over the disk $B(z_0,\varepsilon)$. Note that the
condition (\ref{FMO_eq2.4}) includes the assumption that $\varphi$
is integrable in some neighborhood of the point $z_0$. We say also
that a function $\varphi:D\to{\Bbb R}$ is of {\bf finite mean
oscillation in $D$}, abbr. $\varphi\in{\rm FMO}(D)$ or simply
$\varphi\in{\rm FMO}$, if $\varphi\in{\rm FMO}(z_0)$ for all points
$z_0\in D$. We write $\varphi\in{\rm FMO}(\overline{D})$ if
$\varphi$ is given in a domain $G$ in $\Bbb{C}$ such that
$\overline{D}\subset G$ and $\varphi\in{\rm FMO}(G)$.

\medskip

The following statement is obvious by the triangle inequality.


\begin{propo}\label{FMO_pr2.1} If, for a  collection of numbers
$\varphi_{\varepsilon}\in{\Bbb R}$,
$\varepsilon\in(0,\varepsilon_0]$,
\begin{equation}\label{FMO_eq2.7}\overline{\lim\limits_{\varepsilon\to0}}\ \ \
\dashint_{B(z_0,\varepsilon)}|\varphi(z)-\varphi_{\varepsilon}|\,dm(z)<\infty\,,\end{equation}
then $\varphi $ is of finite mean oscillation at $z_0$.
\end{propo}


In particular choosing here  $\varphi_{\varepsilon}\equiv0$,
$\varepsilon\in(0,\varepsilon_0]$, we obtain the following.


\begin{corollary}\label{FMO_cor2.1} If, for a point $z_0\in D$,
\begin{equation}\label{FMO_eq2.8}\overline{\lim\limits_{\varepsilon\to 0}}\ \ \
\dashint_{B(z_0,\varepsilon)}|\varphi(z)|\,dm(z)<\infty\,,
\end{equation} then $\varphi$ has finite mean oscillation at
$z_0$. \end{corollary}


Recall that a point $z_0\in D$ is called a {\bf Lebesgue point} of a
function $\varphi:D\to{\Bbb R}$ if $\varphi$ is integrable in a
neighborhood of $z_0$ and \begin{equation}\label{FMO_eq2.7a}
\lim\limits_{\varepsilon\to 0}\ \ \ \dashint_{B(z_0,\varepsilon)}
|\varphi(z)-\varphi(z_0)|\,dm(z)=0\,.\end{equation} It is known
that, almost every point in $D$ is a Lebesgue point for every
function $\varphi\in L^1(D)$. Thus we have by Proposition
\ref{FMO_pr2.1} the following corollary.


\begin{corollary}\label{FMO_cor2.7b} Every
locally integrable function $\varphi:D\to{\Bbb R}$ has a finite mean
oscillation at almost every point in $D$.
\end{corollary}


\begin{remark}\label{FMO_rmk2.13a} Note that the function $\varphi(z)=\log\left(1/|z|\right)$
belongs to BMO in the unit disk $\Delta$, see, e.g., \cite{RR}, p.
5, and hence also to FMO. However,
$\widetilde{\varphi}_{\varepsilon}(0)\to\infty$ as
$\varepsilon\to0$, showing that condition (\ref{FMO_eq2.8}) is only
sufficient but not necessary for a function $\varphi$ to be of
finite mean oscillation at $z_0$. Clearly, ${\rm BMO}(D)\subset{\rm
BMO}_{\rm loc}(D)\subset{\rm FMO}(D)$ and as well-known ${\rm
BMO}_{\rm loc}\subset L_{\rm loc}^p$ for all $p\in[1,\infty)$, see,
e.g., \cite{JN} or \cite{RR}. However, FMO is not a subclass of
$L_{\rm loc}^p$ for any $p>1$ but only of $L_{\rm loc}^1$. Thus, the
class FMO is much more wide than ${\rm BMO}_{\rm loc}$.\end{remark}


Versions of the next lemma has been first proved for BMO in
\cite{RSY$_1$}. For FMO, see the papers \cite{IR, RS, RSY$_7$,
RSY$_4$}  and the monographs \cite{GRSY} and \cite{MRSY}.

\medskip

\begin{lemma}\label{lem13.4.2} Let $D$ be a domain in ${\Bbb C}$ and let
$\varphi:D\to{\Bbb R}$ be a  non-negative function  of the class
${\rm FMO}(z_0)$ for some $z_0\in D$. Then
\begin{equation}\label{eq13.4.5}\int\limits_{\varepsilon<|z-z_0|<\varepsilon_0}\frac{\varphi(z)\,dm(z)}
{\left(|z-z_0|\log\frac{1}{|z-z_0|}\right)^2}=O\left(\log\log\frac{1}{\varepsilon}\right)\
\quad\text{as}\quad\varepsilon\to 0\end{equation} for some
$\varepsilon_0\in(0,\delta_0)$ where $\delta_0=\min(e^{-e},d_0)$,
$d_0=\inf\limits_{z\in\partial D}|z-z_0|$. \end{lemma}


\section{Beltrami equations and ring $Q-$homeomorphisms}

The following notion was motivated by the ring definition of Gehring
for quasiconformal mappings, see, e.g., \cite{Ge$_1$}, and it is
closely relevant with the Beltrami equations. Given a domain $D$ in
${\Bbb C}$ and a Lebesgue measurable function $Q:\Bbb
C\to(0,\infty)$, we say that a homeomorphism $f:D\to\overline{\Bbb
C}$ is a {\bf ring $Q$-homeomorphism at a point} $z_0\in
\overline{D}$ if
\begin{equation}\label{eqOS1.8a}M\left(\Delta\left(fC_1,fC_2;fD\right)\right)\leqslant
\int\limits_{A\cap D}Q(z)\cdot\eta^2(|z-z_0|)\,dm(z)\end{equation}
for any ring $A=A(z_0,r_1,r_2)$ and arbitrary continua $C_1$ and
$C_2$ in $D$ that belong to the different components of the
complement of the ring $A$ in $\overline{\Bbb C}$ including $z_0$
and $\infty$, correspondingly, and for any Lebesgue measurable
function $\eta:(r_1,r_2)\to[0,\infty]$ such that
\begin{equation}\label{eqOS1.9}\int\limits_{r_1}^{r_2}\eta(r)\,dr\geqslant1.\end{equation}

The notion was first introduced at inner points of a domain $D$ in
the work \cite{RSY$_3$}. The ring $Q$-homeomorphisms at boundary
points of a domain $D$ have first been considered in the papers
\cite{RSY333} and \cite{RSY$_6$}.

\medskip

By Lemma 2.2 in \cite{RS} or Lemma 7.4 in \cite{MRSY}, we obtain the
following criterion for homeomorphisms in ${\Bbb C}$ to be ring
$Q$-homeomorphisms, see also Theorem A.7 in \cite{MRSY}.

\begin{lemma}\label{prOS2.2}
{\it Let $D$ and $D'$ be bounded domains in ${\Bbb C}$ and $Q:\Bbb
C\to(0,\infty)$ be a measurable function. A homeomorphism $f:D\to
D'$ is a ring $Q$-homeomorphism at $z_0\in\overline{D}$ if and only
if
\begin{equation}\label{eqOS2.1} M\left(\Delta\left(fS_1,fS_2;fD\right)\right)\leqslant
\left(\int\limits_{r_1}^{r_2}
\frac{dr}{||\,Q||(z_0,r)}\right)^{-1}\quad\quad\forall\
r_1\in(0,r_2)\,,\ r_2\in(0,d_0)\end{equation} where
$S_i=S(z_0,r_i)$, $i=1,2$, $d_0 = \sup\limits_{z\in D}\,|\,z-z_0| $
and $||Q||(z_0,r)$ is the $L_1$-norm of $Q$ over $D\cap S(z_0,r)$.}
\end{lemma}

\medskip

By Theorem 4.1 in \cite{RSSY} every homeomorphic $W^{1,1}_{\rm loc}$
solution of the Beltrami equation (\ref{eqBeltrami}) in a domain
$D\subseteq{\Bbb C}$ is the so--called lower $Q$-homeomorphism at
every point $z_0\in\overline{D}$ with $Q(z)=K^T_{\mu}(z,z_0)$, $z\in
D$, and $Q(z)\equiv \varepsilon > 0$ in $\Bbb C\setminus D$. On the
other hand, by Theorem 2 in \cite{KPRS} for a locally integrable
$Q$, if $f:D\to D'$ is a lower $Q$-homeomorphism at a point
$z_0\in\overline{D}$, then $f$ is a ring $Q$-homeomorphism at the
point $z_0$. Thus, we have the following conclusion.

\medskip

\begin{theorem}\label{thKPR3.1} {\it Let $f$ be a homeomorphic $W^{1,1}_{\rm
loc}$ solution of the Beltrami equation (\ref{eqBeltrami}) in a
domain $D\subseteq{\Bbb C}$ and $K_{\mu}\in L^1(D)$. Then $f$ is a
ring $Q$-homeomorphism at every point $z_0\in\overline{D}$ with
$Q(z)=K^T_{\mu}(z,z_0)$, $z\in D$.}
\end{theorem}

In fact, it is sufficient to assume here instead of the condition
$K_{\mu}\in L^1(D)$ that $K^T_{\mu}(z,z_0)$ is integrable along the
circles $|z-z_0|=r$ for a.e. small enough $r$.

\bigskip

\section{Continuous extension of ring $Q-$homeomorphisms}

\begin{lemma}\label{lemOSKRSS1000inverse} {\it Let $D$ and $D'$ be bounded finitely
connected domains in ${\Bbb C}$ and let $f:D\to D'$ be a ring
$Q_{z_0}$-homeomorphism at every point $z_0\in\partial{D}$. Suppose
that
\begin{equation}\label{eqOSKRSS1000inverse}
\int\limits_{D(z_0,\varepsilon, \varepsilon_0)}
Q_{z_0}(z)\cdot\psi_{z_0,\varepsilon, \varepsilon_0}^2(|z-z_0|)\
dm(z) = o\left(I_{z_0, \varepsilon_0}^2(\varepsilon)\right)\ \
\mbox{as}\ \ \varepsilon\to0\ \ \ \forall\ z_0\in\partial D
\end{equation} where $D(z_0,\varepsilon, \varepsilon_0)=\{z\in
D:\varepsilon<|z-z_0|<\varepsilon_0\}$ for every small enough
$0<\varepsilon_0<d(z_0)=\sup\limits_{z\in D}\,|z-z_0|$ and where
$\psi_{z_0,\varepsilon, \varepsilon_0}(t): (0,\infty)\to
[0,\infty]$, $\varepsilon\in(0,\varepsilon_0)$, is a family of
(Lebesgue) measurable functions such that
$$0\ <\ I_{z_0, \varepsilon_0}(\varepsilon)\ :=\
\int\limits_{\varepsilon}^{\varepsilon_0}\psi_{z_0,\varepsilon,
\varepsilon_0}(t)\ dt\ <\ \infty\qquad\qquad\forall\
\varepsilon\in(0,\varepsilon_0)\ .$$ Then $f$ can be extended to a
continuous mapping of $\ \overline{D}_P$ onto $\
\overline{D^{\prime}}_P$.}
\end{lemma}

\medskip

{\bf Proof.} By Remark \ref{METRIC} with no loss of generality we
may assume that $D^{\prime}$ is a circular domain and, thus,
$\overline{D^{\prime}}_P=\overline{D^{\prime}}$. Let us first prove
that the cluster set
$$L\ =\ C(P,f)\ :=\ \left\{\ \zeta\in{{\Bbb
C}}:\ \zeta\ =\ \lim\limits_{n\to\infty}f(z_n),\ z_n\to P,\ z_n\in
D,\ n=1,2,\ldots\ \right\}$$ consists of a single point
$\zeta_0\in\partial D^{\prime}$ for each prime end $P$ of the domain
$D$.

Note that $L\neq\varnothing$ by compactness of the set
$\overline{D^{\prime}}$, and $L$ is a subset of $\partial
D^{\prime}$, see, e.g., Proposition 2.5 in \cite{RSal} or
Proposition 13.5 in \cite{MRSY}. Let us assume that there is at
least two points $\zeta_0$ and $\zeta_*\in L$. Set
$U=B(\zeta_0,\rho_0)=\{ \zeta\in\Bbb C: |\zeta -\zeta_0|<\rho_0\}$
where $0<\rho_0<|\zeta_*-\zeta_0|$.

Let $\sigma_k$, $k=1,2,\ldots\,$, be a chain of cross--cuts of $D$
in the prime end $P$ lying on circles $S_k=S(z_0,r_k)$ from Lemma
\ref{thabc3} where $z_0 \in \partial D.$ Let $D_k$, $k=1,2,\ldots $
be the domains associated with $\sigma_k$. Then there exist points
$\zeta_k$ and $\zeta^*_k$ in the domains $D_{k}'=f(D_{k})$ such that
$|\zeta_0-\zeta_k|<\rho_0$ and $|\zeta_0-\zeta^*_k|>\rho_0$ and,
moreover, $\zeta_k\to \zeta_0$ and $\zeta^*_k\to \zeta_*$ as
$k\to\infty$. Let $C_k$ be continuous curves joining $\zeta_k$ and
$\zeta^*_k$ in $D_{k}'$. Note that by the construction $\partial
U\cap C_k\neq\varnothing$, $k=1,2,\ldots $.

By the condition of strong accessibility of the point $\zeta_0$,
there is a continuum $E\subset D'$ and a number $\delta>0$ such that
\begin{equation}\label{delta} M(\Delta(E,C_k;D'))\ \geqslant\ \delta\end{equation}
for all large enough $k$. Note that $C=f^{-1}(E)$ is a compact
subset of $D$ and hence $d_0={\rm dist}(z_0,C)>0$. Let
$\varepsilon_0\in(0,d_0)$ be small enough from the hypotheses of the
lemma. With no loss of generality, we may assume that
$r_k<\varepsilon_0$ and (\ref{delta}) holds for all $k=1,2,\ldots$.

Let $\Gamma_{m}$ be a family of all continuous curves in $D\setminus
D_m$ joining the circle $S_{0}=S(z_0,\varepsilon_0)$ and
$\overline{\sigma_m}$, $m=1,2,\ldots$. Note that by the construction
$C_k\subset D_k^{\prime}\subset D_m'$ for all $k\geqslant m$ and,
thus, by the principle of minorization
$M(f(\Gamma_{m}))\geqslant\delta$ for all $m=1,2,\ldots$.

On the other hand, every function $\eta(t)=\eta_m(t):=
\psi_{z_0,r_m,\varepsilon_0}^*(t)/I_{z_0,\varepsilon_0}(r_m)$,
$m=1,2,\ldots$, satisfies the condition (\ref{eqOS1.9}) and hence
$$M(f\Gamma_m)\ \leqslant\ \int\limits_{D}Q(z)\cdot
{\eta^2_m}(z)\ dm(z)\ ,$$ i.e., $M(f\Gamma_m)\to 0$ as $m\to\infty$
in view of (\ref{eqOSKRSS1000inverse}).


The obtained contradiction disproves the assumption that the cluster
set $C(P,f)$ consists of more than one point.


Thus, we have the extension $h$ of $f$ to $\overline{D}_P$ such that
$C(E_D, f)\subseteq \partial{D^{\prime}}$. In fact, $C(E_D, f)=
\partial{D^{\prime}}$. Indeed, if $\zeta_0\in D^{\prime}$,
then there is a sequence $\zeta_n$ in $D^{\prime}$ being convergent
to $\zeta_0$. We may assume with no loss of generality that
$f^{-1}(\zeta_n)\to P_0\in E_D$ because $\overline{D}_P$ is compact,
see Remark \ref{METRIC}. Hence $\zeta_0\in C(P_0, f)$.

Finally, let us show that the extended mapping
$h:\overline{D}_P\to\overline{D^{\prime}}$ is continuous. Indeed,
let $P_n\to P_0$ in $\overline{D}_P$. If $P_0\in D$, then the
statement is obvious. If $P_0\in E_D$, then by the last item we are
able to choose $P^*_n\in D$ such that $\rho(P_n,P^*_n)<1/n$ where
$\rho$ is one of the metrics in Remark \ref{METRIC} and
$|h(P_n)-h(P^*_n)|<1/n$. Note that by the first part of the proof
$h(P^*_n)\to h(P_0)$ because $P^*_n\to P_0$. Consequently,
$h(P_n)\to h(P_0)$, too. $\Box$

\begin{theorem}\label{l:6.3c} {\it
Let $D$ and $D'$ be bounded finitely connected domains in ${\Bbb C}$
and let $f:D\to D'$ be a ring $Q_{z_0}$-homeomorphism at every point
$z_0\in\partial{D}$. If
\begin{equation}\label{e:6.4c}\int\limits_{0}^{\varepsilon(z_0)}
\frac{dr}{||Q_{z_0}||(r)}\ =\ \infty\qquad\qquad \forall\
z_0\in\partial D\end{equation} where
$0<\varepsilon(z_0)<d(z_0):=\sup\limits_{z\in D}\,|z-z_0|$ and
\begin{equation}\label{e:6.6c} ||Q_{z_0}||(r)\ :=\int\limits_{
D\cap S({z_0},r)}Q_{z_0}\ ds\ ,\end{equation} then $f$ can be
extended to a continuous mapping of $\overline{D}_P$ onto
$\overline{D^{\prime}}_P$.}
\end{theorem}

\medskip

{\bf Proof.} Indeed, condition (\ref{e:6.4c}) implies that
\begin{equation}\label{e:6.4dc}\int\limits_{0}^{\varepsilon_0}
\frac{dr}{||Q_{z_0}||(r)}\ =\ \infty\qquad\qquad \forall\
z_0\in\partial D \qquad\forall\
\varepsilon_0\in(0,\varepsilon(z_0))\end{equation} because the left
hand side in (\ref{eqOS2.1}) is not equal to zero, see Theorem 5.2
in \cite{Na$_2$}, and hence by Lemma \ref{prOS2.2}
$$\int\limits_{\varepsilon_0}^{\varepsilon(z_0)}
\frac{dr}{||Q_{z_0}||(r)}\ <\ \infty\ .$$

\medskip

On the other hand, for the functions
\begin{equation}\label{3.21}\psi_{z_0, \varepsilon_0}(t)\ :=\ \left \{\begin{array}{lr}
1/||Q_{z_0}||(t),\quad & \ t\in (0,\varepsilon_0),
\\ 0,  & \ t\in [\varepsilon_0,\infty),\end{array}\right.\end{equation} we
have by the Fubini theorem that
\begin{equation}\label{3.22}\int\limits_{D(z_0,\varepsilon, \varepsilon_0)} Q_{z_0}(z)\cdot\psi_{z_0,\varepsilon_0}^2(|z-z_0|)\
dm(z)\ =\ \int\limits_{\varepsilon}^{\varepsilon_0}
\frac{dr}{||Q_{z_0}||(r)}
\end{equation} and, consequently, condition (\ref{eqOSKRSS1000inverse}) holds by (\ref{e:6.4dc}) for all $z_0\in\partial D$ and
$\varepsilon_0\in(0,\varepsilon(z_0))$.

\medskip

Here we have used the standard conventions in the integral theory
that $a/\infty=0$ for $a\neq\infty$ and $0\cdot\infty=0$, see, e.g.,
Section I.3 in \cite{Sa}.

\bigskip

Thus, Theorem \ref{l:6.3c} follows immediately from Lemma
\ref{lemOSKRSS1000inverse}. $\Box$

\bigskip

\section{Extension of the inverse mappings to the boundary}

The base for the proof on extending the inverse mappings of ring
$Q$-ho\-meo\-mor\-phism by prime ends in the plane is the following
fact on the cluster sets.


\begin{lemma}\label{l:9.1} {\it Let $D$ and $D'$ be bounded finitely
connected domains in ${\Bbb C}$, $P_0$ and $P_*$ be prime ends of
$D$,  $P_*\ne P_0$. Denote  by $\sigma_m$, $m=1,2,\ldots$, a chain
of cross--cuts in $P_0$ from Lemma \ref{thabc3} lying on circles
$S(z_0,r_m)$, $z_0\in\partial D$, with associated domains $d_m$.
Suppose that $Q$ is integrable over $D\cap S(z_0,r)$ for a set $E$
of numbers $r\in(0,\delta)$ of a positive linear measure where
$\delta=r_{m_0}$ and $m_0$ is such that the domain $d_{m_0}$ does
not contain sequences of points converging to $P_*$. If $f:D\to D'$
is a ring $Q$-homeomorphism at the point $z_0$ and $\partial D'$ is
weakly flat, then\begin{equation}\label{e:9.2} C(P_0,f)\cap
C(P_*,f)\ =\ \varnothing\ .\end{equation}}
\end{lemma}

Note that in view of metrizability of the completion
$\overline{D}_P$ of the domain $D$ with prime ends, see Remark
\ref{METRIC}, the number $m_0$ in Lemma \ref{l:9.1} always exists.

\medskip

{\bf Proof.} Let us choose $\varepsilon\in(0,\delta)$ such that
$E_0:=\{r\in E:r\in(\varepsilon,\delta)\}$ has a positive linear
measure. Such a choice is possible in view of subadditivity of the
linear measure and the exhaustion $E = \cup E_n$ where $E_n = \{r\in
E:r\in(1/n,\delta)\}\,,$ $n=1,2,\ldots $. Note that by Lemma
\ref{prOS2.2}, for $S_1=S(z_0,\varepsilon)$ and $S_2=S(z_0,\delta)$,
\begin{equation}\label{e:9.3}
M\left(\Delta\left(fS_1,fS_2;fD\right)\right)\ <\ \infty\ .
\end{equation}

Let us assume that $C_0\cap C_*\neq\varnothing$ where $C_0=C(P_0,f)$
and $C_*=C(P_*,f)$. By the construction there is $m_1>m_0$ such that
$\sigma_{m_1}$ lies on the circle $S(z_0,r_{m_1})$ with
$r_{m_1}<\varepsilon$. Let $d_0=d_{m_1}$ and $d_*\subseteq
D\setminus d_{m_0}$ be a domain associated with a chain of
cross--cuts in the prime end  $P_*$. Let $\zeta_0\in C_1\cap C_2$.
Choose $\rho_0>0$ such that $S(\zeta_0,\rho_0)\cap
f(d_0)\neq\varnothing$ and $S(\zeta_0,\rho_0)\cap
f(d_*)\neq\varnothing$.

Set $\Gamma=\Delta(\overline{d_0},\overline{d_*};D)$.
Correspondingly (\ref{e:9.3}), by the principle of minorization
\begin{equation}\label{e:9.4}M(f(\Gamma))\  <\ \infty\ .
\end{equation}
Let $M_0>M(f(\Gamma))$ be a finite number. By the condition of the
lemma, $\partial D'$ is weakly flat and hence there is
$\rho_*\in(0,\rho_0)$ such that
$$M(\Delta(E,F;D'))\ \geqslant\ M_0$$
for all continua $E$ and $F$ in $D'$ intersecting the circles
$S(\zeta_0,\rho_0)$ and $S(\zeta_0,\rho_*)$. However, these circles
can be joined by continuous curves $c_1$ and $c_2$ in the domains
$f(d_0)$ and $f(d_*)$, correspondingly, and, in particular, for
these curves
\begin{equation}\label{e:9.4a}M_0\ \leqslant\
M(\Delta(c_1,c_2;D'))\ \leqslant\ M(f(\Gamma))\,.\end{equation} The
obtained contradiction disproves the assumption that $C_0\cap
C_*\neq\varnothing$. $\Box$

\medskip

\begin{theorem}\label{thKPR8.2} {\it Let $D$ and $D'$ be bounded finitely
connected domains in ${\Bbb C}$ and $f:D\to D'$ be a
$Q_{z_0}$-homeomorphism at every point $z_0\in\partial{D}$ with
$Q_{z_0}\in L^{1}(D\cap U_{z_0})$ for a neighborhood $U_{z_0}$ of
$z_0$. Then $f^{-1}$ can be extended to a continuous mapping of
$\overline{D^{\prime}}_P$ onto $\overline{D}_P$.}
\end{theorem}

\medskip

{\it Proof.} By Remark \ref{METRIC}, we may assume with no loss of
generality that $D^{\prime}$ is a circular domain,
$\overline{D^{\prime}}_P=\overline{D^{\prime}}$; $C(\zeta_0,
f^{-1})\ne\varnothing $ for every $\zeta_0\in
\partial D^{\prime}$ because $\overline{D}_P$ is metrizable and compact.
Moreover, $C(\zeta_0, f^{-1})\cap D=\varnothing $, see, e.g.,
Proposition 2.5 in \cite{RSal} or Proposition 13.5 in \cite{MRSY}.

Let us assume that there is at least two different prime ends $P_1$
and $P_2$ in $C(\zeta_0, f^{-1})$. Then $\zeta_0\in C(P_1,f)\cap
C(P_2,f)$. Let $z_1\in\partial D$ be a point corresponding to $P_1$
from Lemma \ref{thabc3}. Note that
\begin{equation}\label{eqKPR6.2ad} E\ =\ \{r\in(0,\delta):\  Q_{z_1}|_{D\cap S(z_1,r)}\in L^{1}(D\cap S(z_1,r))\}
\end{equation} has a positive linear measure for every $\delta>0$ by the Fubini theorem, see, e.g., \cite{Sa}, because $Q_{z_1}\in
L^{1}(D\cap U_{z_1})$.  The obtained contradiction with Lemma
\ref{l:9.1} shows that $C(\zeta_0, f^{-1})$ contains only one prime
end of $D$.

Thus, we have the extension $g$ of $f^{-1}$ to
$\overline{D^{\prime}}$ such that $C(\partial D^{\prime},
f^{-1})\subseteq \overline{D}_P\setminus D$. In fact, $C(\partial
D^{\prime}, f^{-1})=\overline{D}_P\setminus D$. Indeed, if $P_0$ is
a prime end of $D$, then there is a sequence $z_n$ in $D$ being
convergent to $P_0$. We may assume without loss of generality that
$z_n\to z_0\in\partial D$ and $f(z_n)\to \zeta_0\in\partial
D^{\prime}$ because $\overline{D}$ and $\overline{D^{\prime}}$ are
compact. Hence $P_0\in C(\zeta_0, f^{-1})$.

Finally, let us show that the extended mapping
$g:\overline{D^{\prime}}\to\overline{D}_P$ is continuous. Indeed,
let $\zeta_n\to\zeta_0$ in $\overline{D^{\prime}}$. If $\zeta_0\in
D^{\prime}$, then the statement is obvious. If $\zeta_0\in\partial
D^{\prime}$, then take $\zeta^*_n\in D^{\prime}$ such that
$|\zeta_n-\zeta^*_n|<1/n$ and $\rho(g(\zeta_n),g(\zeta^*_n))<1/n$
where $\rho$ is one of the metrics in Remark \ref{METRIC}. Note that
by the construction $g(\zeta^*_n)\to g(\zeta_0)$ because
$\zeta^*_n\to \zeta_0$. Consequently, $g(\zeta_n)\to g(\zeta_0)$,
too. $\Box$

\medskip

\begin{theorem}\label{t:9.12} {\it Let $D$ and $D'$ be bounded finitely
connected domains in ${\Bbb C}$ and $f:D\to D'$ be a
$Q_{z_0}$-homeomorphism at every point $z_0\in\partial{D}$ with the
condition
\begin{equation}\label{e:6.4}\int\limits_{0}^{\varepsilon(z_0)}
\frac{dr}{||Q_{z_0}||_{}(r)}\ =\ \infty\end{equation} where
$0<\varepsilon(z_0)<d(z_0)=\sup\limits_{z\in D}\,|z-z_0|$ and
\begin{equation}\label{e:6.6} ||Q_{z_0}||(r)\ =\int\limits_{
D\cap S({z_0},r)}Q_{z_0}\ ds\ .\end{equation} Then $f^{-1}$ can be
extended to a continuous mapping of $\overline{D^{\prime}}_P$ onto
$\overline{D}_P$.}\end{theorem}

\medskip

{\bf Proof.} Indeed, condition (\ref{e:6.4}) implies that
\begin{equation}\label{e:6.4d}\int\limits_{0}^{\delta}
\frac{dr}{||Q_{z_0}||(r)}\ =\ \infty\qquad\qquad \forall\
z_0\in\partial D \qquad\forall\
\delta\in(0,\varepsilon(z_0))\end{equation} because the left hand
side in (\ref{eqOS2.1}) is not equal to zero, see Theorem 5.2 in
\cite{Na$_2$}, and hence by Lemma \ref{prOS2.2}
$$\int\limits_{\delta}^{\varepsilon(z_0)}
\frac{dr}{||Q_{z_0}||(r)}\ <\ \infty\ .$$ Thus, the set
\begin{equation}\label{eqKPR6.2ad} E\ =\ \{r\in(0,\delta):\  Q_{z_0}|_{D\cap S(z_0,r)}\in L^{1}(D\cap S(z_0,r))\}
\end{equation}
has a positive linear measure for all $z_0\in\partial D$ and all
$\delta\in(0,\varepsilon(z_0))$. The rest of arguments is perfectly
similar to one in the proof of Theorem \ref{thKPR8.2}. $\Box$

\bigskip

\section{Homeomorphic extension of ring $Q-$homeomorphisms}

Combining Theorems \ref{l:6.3c} and \ref{t:9.12}, we obtain the
following principal result.


\begin{theorem}\label{t:10.1} {\it Let $D$ and $D'$ be bounded finitely connected domains in ${\Bbb C}$
and let $f:D\to D'$ be a ring $Q_{z_0}$-homeomorphism at every point
$z_0\in\partial{D}$. If
\begin{equation}\label{e:6.4ch}\int\limits_{0}^{\varepsilon(z_0)}
\frac{dr}{||Q_{z_0}||(r)}\ =\ \infty\qquad\qquad \forall\
z_0\in\partial D\end{equation} where
$0<\varepsilon(z_0)<d(z_0):=\sup\limits_{z\in D}\,|z-z_0|$ and
\begin{equation}\label{e:6.6ch} ||Q_{z_0}||(r)\ :=\int\limits_{
D\cap S({z_0},r)}Q_{z_0}\ ds\ ,\end{equation} then $f$ can be
extended to a homeomorphism of $\overline{D}_P$ onto
$\overline{D^{\prime}}_P$.}
\end{theorem}


\begin{corollary}\label{thOSKRSS100} {\it In particular, the conclusion of Theorem
\ref{t:10.1} holds if
\begin{equation}\label{eqOSKRSS100d}q_{z_0}(r)=O\left(\log{\frac1r}\right)\qquad\qquad  \forall\
z_0\in\partial D \end{equation} as $r\to0$ where $q_{z_0}(r)$ is the
average of $Q$ over the circle $|z-z_0|=r$. }
\end{corollary}


Using Lemma 2.2 in \cite{RS}, see also Lemma 7.4 in \cite{MRSY}, by
Theorem \ref{t:10.1} we obtain the following general lemma that, in
turn, makes possible to obtain new criteria in a great number.


\begin{lemma}\label{lemOSKRSS1000} {\it Let $D$ and $D'$ be bounded finitely
connected domains in ${\Bbb C}$ and let $f:D\to D'$ be a ring
$Q_{z_0}$-homeomorphism at every point $z_0\in\partial{D}$ where
$Q_{z_0}$ is integrable in a neighborhood of $z_0$. Suppose that
\begin{equation}\label{eqOSKRSS1000inverseh}
\int\limits_{D(z_0,\varepsilon, \varepsilon_0)}
Q_{z_0}(z)\cdot\psi_{z_0,\varepsilon, \varepsilon_0}^2(|z-z_0|)\
dm(z) = o\left(I_{z_0, \varepsilon_0}^2(\varepsilon)\right)\ \
\mbox{as}\ \ \varepsilon\to0\ \ \ \forall\ z_0\in\partial D
\end{equation} where $D(z_0,\varepsilon, \varepsilon_0)=\{z\in
D:\varepsilon<|z-z_0|<\varepsilon_0\}$ for every small enough
$0<\varepsilon_0<d(z_0)=\sup\limits_{z\in D}\,|z-z_0|$ and where
$\psi_{z_0,\varepsilon, \varepsilon_0}(t): (0,\infty)\to
[0,\infty]$, $\varepsilon\in(0,\varepsilon_0)$, is a family of
(Lebesgue) measurable functions such that
$$0\ <\ I_{z_0, \varepsilon_0}(\varepsilon)\ :=\
\int\limits_{\varepsilon}^{\varepsilon_0}\psi_{z_0,\varepsilon,
\varepsilon_0}(t)\ dt\ <\ \infty\qquad\qquad\forall\
\varepsilon\in(0,\varepsilon_0)\ .$$ Then $f$ can be extended to a
 homeomorphism of $\ \overline{D}_P$ onto $\
\overline{D^{\prime}}_P$.}
\end{lemma}


\begin{remark}\label{rmKR2.9} In fact, instead of
integrability of $Q_{z_0}$ in a neighborhood of $z_0$, it is
sufficient to request that $Q_{z_0}$ is integrable over $D\cap
S(z_0,r)$ for a.e. $r\in(0,\varepsilon_0)$.

Note also that it is not only Lemma \ref{lemOSKRSS1000} follows from
Theorem \ref{t:10.1} under the given conditions on integrability of
$Q_{z_0}$ but, inversely, Theorem \ref{t:10.1} follows from Lemma
\ref{lemOSKRSS1000}, too, as it was shown under the proof of Theorem
\ref{l:6.3c}. Thus, Theorem \ref{t:10.1} is equivalent under the
given conditions to Lemma \ref{lemOSKRSS1000} but each of them is
sometimes more convenient for applications than another one.

Finally, note that (\ref{eqOSKRSS1000inverseh}) holds, in
particular, if
\begin{equation}\label{eqOSKRSS100a} \int\limits_{D(z_0,
\varepsilon_0)}Q_{z_0}(z)\cdot\psi^2 (|z-z_0|)\ dm(z)\ <\
\infty\qquad\qquad \forall\ z_0\in\partial D\end{equation} where
$D(z_0,\varepsilon_0)=\{z\in D:|z-z_0|<\varepsilon_0\}$ and where
$\psi(t): (0,\infty)\to [0,\infty]$ is a locally integrable function
such that $I_{z_0,\varepsilon_0}(\varepsilon)\to\infty$ as
$\varepsilon\to0$. In other words, for the extendability of $f$ to a
homeomorphism of $\overline{D}_P$ onto $\overline{D'}_P$, it
suffices for the integrals in (\ref{eqOSKRSS100a}) to be convergent
for some nonnegative function $\psi(t)$ that is locally integrable
on $(0,\infty)$ but that has a non-integrable singularity at zero.
\end{remark}


Choosing in Lemma \ref{lemOSKRSS1000} $\psi(t):=\frac{1}{t\log 1/t}$
and applying Lemma \ref{lem13.4.2}, we obtain the next result.


\begin{theorem}\label{thOSKRSS101} {\it Let $D$ and $D'$ be bounded finitely connected domains in ${\Bbb C}$
and let $f:D\to D'$ be a ring $Q_{z_0}$-homeomorphism at every point
$z_0\in\partial{D}$ where $Q_{z_0}$ has finite mean oscillation at
$z_0$. Then $f$ can be extended to a homeomorphism of
$\overline{D}_P$ onto $\overline{D'}_P$.}
\end{theorem}


\begin{corollary}\label{corOSKRSS6.6.2} {\it In particular, the conclusion of Theorem \ref{thOSKRSS101} holds if
\begin{equation}\label{eqOSKRSS6.6.3}
\overline{\lim\limits_{\varepsilon\to0}}\ \
\dashint_{B(z_0,\varepsilon)}Q_{z_0}(z)\ dm(z)\ <\
\infty\qquad\qquad \forall\ z_0\in\partial D\end{equation}}
\end{corollary}


\begin{corollary}\label{corOSKRSS6.6.33} {\it The conslusion of Theorem \ref{thOSKRSS101} holds if
every point $z_0\in\partial D$ is a Lebesgue point of the function
$Q_{z_0}$.}
\end{corollary}


The next statement also follows from Lemma \ref{lemOSKRSS1000} under
the choice $\psi(t)=1/t.$


\begin{theorem}\label{thOSKRSS102} {\it Let $D$ and $D'$ be
bounded finitely connected domains in ${\Bbb C}$ and $f:D\to D'$ be
a ring $Q_{z_0}$-homeomorphism at every point $z_0\in\partial{D}$.
If, for some $\varepsilon_0=\varepsilon(z_0)>0$,
\begin{equation}\label{eqOSKRSS10.336a}\int\limits_{\varepsilon<|z-z_0|<\varepsilon_0}Q_{z_0}(z)\
\frac{dm(z)}{|z-z_0|^2}\ =\
o\left(\left[\log\frac{1}{\varepsilon}\right]^2\right)\qquad
\mbox{as}\ \ \ \varepsilon\to 0\qquad \forall\ z_0\in\partial
D\end{equation} then $f$ can be extended to a homeomorphism of $\
\overline{D}_P$ onto $\overline{D'}_P$.}
\end{theorem}


\begin{remark}\label{rmOSKRSS200} Choosing in Lemma \ref{lemOSKRSS1000}
the function $\psi(t)=1/(t\log 1/t)$ instead of $\psi(t)=1/t$,
(\ref{eqOSKRSS10.336a}) can be replaced by the more weak condition
\begin{equation}\label{eqOSKRSS10.336b}
\int\limits_{\varepsilon<|z-z_0|<\varepsilon_0}\frac{Q_{z_0}(z)\
dm(z)}{\left( |z-z_0|\ \log{\frac{1}{|z-z_0|}}\right) ^2}\ =\
o\left(\left[\log\log\frac{1}{\varepsilon}\right]^2\right)\end{equation}
and (\ref{eqOSKRSS100d}) by the condition
\begin{equation}\label{eqOSKRSS10.336h} q_{z_0}(r)\ =\ o
\left(\log\frac{1}{r}\log\,\log\frac{1}{r} \right).\end{equation} Of
course, we could give here the whole scale of the corresponding
condition of the logarithmic type using suitable functions
$\psi(t).$
\end{remark}


Theorem \ref{t:10.1} has a magnitude of other fine consequences, for
instance:


\begin{theorem}\label{thOSKRSS103} {\it Let $D$ and $D'$ be bounded finitely
connected domains in ${\Bbb C}$ and let $f:D\to D'$ be a  ring
$Q_{z_0}$-homeomorphism at every point $z_0\in\partial{D}$ and
\begin{equation}\label{eqOSKRSS10.36b} \int\limits_{D\cap B(z_0,\varepsilon_0)}\Phi_{z_0}\left(Q_{z_0}(z)\right)\ dm(z)\ <\ \infty
\qquad \forall\ z_0\in\partial D\end{equation} for
$\varepsilon_0=\varepsilon(z_0)>0$ and a nondecreasing convex
function $\Phi_{z_0}:[0,\infty)\to[0,\infty)$ with
\begin{equation}\label{eqOSKRSS10.37b}
\int\limits_{\delta(z_0)}^{\infty}\frac{d\tau}{\tau\Phi_{z_0}^{-1}(\tau)}\
=\ \infty\end{equation} for $\delta(z_0)>\Phi_{z_0}(0)$. Then $f$ is
extended to a homeomorphism of $\overline{D}_P$ onto
$\overline{D'}_P$.}
\end{theorem}


Indeed, by Theorem 3.1 and Corollary 3.2 in
 \cite{RSY}, (\ref{eqOSKRSS10.36b}) and
(\ref{eqOSKRSS10.37b}) imply (\ref{e:6.4ch}) and, thus, Theorem
\ref{thOSKRSS103} is a direct consequence of Theorem \ref{t:10.1}.


\begin{corollary}\label{corOSKRSS6.6.3} {\it In particular, the conclusion of Theorem
\ref{thOSKRSS103} holds if
\begin{equation}\label{eqOSKRSS6.6.6}
\int\limits_{D\cap B(z_0,\varepsilon_0)}e^{\alpha_0 Q_{z_0}(z)}\
dm(z)\ <\ \infty\qquad \forall\ z_0\in\partial D\end{equation} for
some $\varepsilon_0=\varepsilon(z_0)>0$ and
$\alpha_0=\alpha(z_0)>0$.}
\end{corollary}


\begin{remark}\label{rmOSKRSS200000}
By Theorem 2.1 in \cite{RSY}, see also Proposition 2.3 in
\cite{RS1}, (\ref{eqOSKRSS10.37b}) is equivalent to every of the
conditions from the following series:
\begin{equation}\label{eq333Y}\int\limits_{\delta(z_0)}^{\infty}
H'_{z_0}(t)\ \frac{dt}{t}=\infty\ ,\quad\ \delta(z_0)>0\
,\end{equation}
\begin{equation}\label{eq333F}\int\limits_{\delta(z_0)}^{\infty}
\frac{dH_{z_0}(t)}{t}=\infty\ ,\quad\ \delta(z_0)>0\ ,\end{equation}
\begin{equation}\label{eq333B}
\int\limits_{\delta(z_0)}^{\infty}H_{z_0}(t)\ \frac{dt}{t^2}=\infty\
,\quad\ \delta(z_0)>0\ ,
\end{equation}
\begin{equation}\label{eq333C}
\int\limits_{0}^{\Delta(z_0)}H_{z_0}\left(\frac{1}{t}\right)\,dt=\infty\
,\quad\ \Delta(z_0)>0\ ,
\end{equation}
\begin{equation}\label{eq333D}
\int\limits_{\delta_*(z_0)}^{\infty}\frac{d\eta}{H_{z_0}^{-1}(\eta)}=\infty\
,\quad\ \delta_*(z_0)>H_{z_0}(0)\ ,
\end{equation}
where
\begin{equation}\label{eq333E}
H_{z_0}(t)=\log\Phi_{z_0}(t)\ .\end{equation}

Here the integral in (\ref{eq333F}) is understood as the
Lebesgue--Stieltjes integral and the integrals in (\ref{eq333Y}) and
(\ref{eq333B})--(\ref{eq333D}) as the ordinary Lebesgue integrals.


It is necessary to give one more explanation. From the right hand
sides in the conditions (\ref{eq333Y})--(\ref{eq333D}) we have in
mind $+\infty$. If $\Phi_{z_0}(t)=0$ for $t\in[0,t_*(z_0)]$, then
$H_{z_0}(t)=-\infty$ for $t\in[0,t_*(z_0)]$ and we complete the
definition $H'_{z_0}(t)=0$ for $t\in[0,t_*(z_0)]$. Note, the
conditions (\ref{eq333F}) and (\ref{eq333B}) exclude that $t_*(z_0)$
belongs to the interval of integrability because in the contrary
case the left hand sides in (\ref{eq333F}) and (\ref{eq333B}) are
either equal to $-\infty$ or indeterminate. Hence we may assume in
(\ref{eq333Y})--(\ref{eq333C}) that $\delta(z_0)>t_0$,
correspondingly, $\Delta(z_0)<1/t(z_0)$ where
$t(z_0):=\sup\limits_{\Phi_{z_0}(t)=0}t$, set $t(z_0)=0$ if
$\Phi_{z_0}(0)>0$.


The most interesting of the above conditions is (\ref{eq333B}) that
can be rewritten in the form:
\begin{equation}\label{eq5!}
\int\limits_{\delta(z_0)}^{\infty}\log \Phi_{z_0}(t)\ \
\frac{dt}{t^{2}}\ =\ \infty\ .
\end{equation}

Note also that, under every homeomorphism $f$ between domains $D$
and $D^{\prime}$ in $\overline{{\Bbb C}}$, there is a natural
one-to-one correspondence between components of their boundaries
$\partial D$ and $\partial D'$, see, e.g., Lemma 5.3 in \cite{IR} or
Lemma 6.5 in \cite{MRSY}. Thus, if a bounded domain $D$ in ${\Bbb
C}$ is finitely connected and $D^{\prime}$ is bounded, then
$D^{\prime}$ is finitely connected, too.

\bigskip

Finally, note that if a domain $D$ in ${\Bbb C}$ is locally
connected on its boundary, then there is a natural one-to-one
correspondence between prime ends of $D$ and boundary points of $D$.
Thus, if $D$ and $D^{\prime}$ are in addition locally connected on
their boundaries in theorems of Sections 7, then $f$ is extended to
a homeomorphism of $\overline D$ onto $\overline{D^{\prime}}$. We
obtained before it similar results when $\partial D^{\prime}$ was
weakly flat which is a more strong condition than local connectivity
of $D^{\prime}$ on its boundary, see, e.g., \cite{KPR} and
\cite{KPRS}.

\bigskip

As known, every Jordan domain $D$ in ${\Bbb C}$ is locally connected
on its boundary, see, e.g., \cite{Wi}, p. 66. It is easy to see, the
latter implies that every bounded finitely connected domain $D$ in
${\Bbb C}$ whose boundary consists of mutually disjoint Jordan
curves and isolated points is also locally connected on its
boundary.

\medskip

Inversely, every bounded finitely connected domain $D$ in ${\Bbb C}$
which is locally connected o
n its boundary has a boundary consisting
of mutually disjoint Jordan curves and isolated points. Indeed,
every such a domain $D$ can be mapped by a conformal mapping $f$
onto the so-called circular domain $D_*$ bounded by a finite
collection of mutually disjoint circles and isolated points, see,
e.g., Theorem V.6.2 in \cite{Goluzin}, that is extended to a
homeomorphism of $\overline D$ onto $\overline{D_*}$, see Remark
\ref{METRIC}.
\end{remark}

\medskip

\section{Boundary behavior of homeomorphic solutions}

On the basis of results in Sections 6 and 7, we obtain by Theorem
\ref{thKPR3.1} the cor\-res\-pon\-ding results on the boundary
behavior of solutions of the Beltrami equations.

\medskip

\begin{theorem}\label{thKPR8.2F} {\it Let $D$ and $D^{\prime}$ be
bounded finitely connected domains in ${\Bbb C}$ and $f:D\to
D^{\prime}$ be a homeomorphic solution in the class $W^{1,1}_{\rm
loc}$ of (\ref{eqBeltrami}) with $K^T_{\mu}(\cdot ,z_0)\in
L^{1}(D\cap B(z_0,\varepsilon_0))$ for every $z_0\in\partial D$.
Then $f^{-1}$ is extended to a continuous mapping of
$\overline{D^{\prime}}_P$ onto $\overline{D}_P$.}
\end{theorem}

\medskip

However, any degree of integrability of $K_{\mu}$ does not guarantee
a continuous extendability of the direct mapping $f$ to the
boundary, see, e.g., an example in the proof of Proposition 6.3 in
\cite{MRSY}. Conditions for it have perfectly another nature. The
principal relevant result is the following.

\medskip

\begin{theorem}\label{t:10.1F} {\it Let $D$ and $D^{\prime}$ be
bounded finitely connected domains in ${\Bbb C}$ and $f:D\to
D^{\prime}$ be a homeomorphic solution in the class $W^{1,1}_{\rm
loc}$ of the Beltrami equation (\ref{eqBeltrami}) with the condition
\begin{equation}\label{ee:6.4F}\int\limits_{0}^{\varepsilon_0}
\frac{dr}{||K^T_{\mu}||(z_0,r)}\ =\ \infty\qquad\forall\
z_0\in\partial D\end{equation} where
$0<\varepsilon_0=\varepsilon(z_0)<d(z_0):=\sup\limits_{z\in
D}\,|z-z_0|$ and
\begin{equation}\label{eq8.7.6F} ||K^T_{\mu}||(z_0,r)\ =\ \int\limits_{
S({z_0},r)}K^T_{\mu}(z,z_0)\ ds\ .\end{equation} Then $f$ can be
extended to a homeomorphism of $\overline{D}_P$ onto
$\overline{D^{\prime}}_P$.}
\end{theorem}

\medskip

Here and later on, we set that $K^T_{\mu}$ is equal to zero outside
of the domain $D$.

\medskip

\begin{corollary}\label{thOSKRSS100F} {\it In particular, the
conclusion of Theorem \ref{t:10.1F} holds if
\begin{equation}\label{eqOSKRSS100dF}k^T_{z_0}(r)=O\left(\log{\frac1r}\right)\
\qquad\forall\ z_0\in\partial D \end{equation} as $r\to0$ where
$k^T_{z_0}(r)$ is the average of $K^T_{\mu}(z,z_0)$ over the circle
$|z-z_0|=r$.}
\end{corollary}

\medskip

\begin{lemma}\label{lemOSKRSS1000F} {\it Let $D$ and $D^{\prime}$ be
bounded finitely connected domains in ${\Bbb C}$ and $f:D\to
D^{\prime}$ be a homeomorphic solution in the class $W^{1,1}_{\rm
loc}$ of the Beltrami equation (\ref{eqBeltrami}) with $K_{\mu}\in
L^{1}(D)$ and
\begin{equation}\label{eqOSKRSS1000F}
\int\limits_{\varepsilon<|z-z_0|<\varepsilon_0}K^T_{\mu}(z,z_0)\cdot\psi_{z_0,\varepsilon}^2(|z-z_0|)\
dm(z)\ =\ o\left(I_{z_0}^2(\varepsilon)\right)\qquad\forall\
z_0\in\partial D
\end{equation} as $\varepsilon\to0$ where
$0<\varepsilon_0<\sup\limits_{z\in D}\,|z-z_0|$ and
$\psi_{z_0,\varepsilon}(t):(0,\infty)\to[0,\infty]$,
$\varepsilon\in(0,\varepsilon_0)$, is a two-parametric family of
measurable functions such that
$$0<I_{z_0}(\varepsilon)\ :=\ \int\limits_{\varepsilon}^{\varepsilon_0}\psi_{z_0,\varepsilon}(t)\ dt\ <\ \infty
\qquad\forall\ \varepsilon\in(0,\varepsilon_0)\ .$$ Then $f$ can be
extended to a homeomorphism of $\overline{D}_P$ onto
$\overline{D^{\prime}}_P$.} \end{lemma}

\medskip

\begin{theorem}\label{thOSKRSS101F} {\it Let $D$ and $D^{\prime}$ be
bounded finitely connected domains in ${\Bbb C}$ and $f:D\to
D^{\prime}$ be a homeomorphic solution in the class $W^{1,1}_{\rm
loc}$ of the Beltrami equation (\ref{eqBeltrami}) with
$K^T_{\mu}(z,z_0)$ of finite mean oscillation at every point
$z_0\in\partial D$. Then $f$ can be extended to a homeomorphism of
$\overline{D}_P$ onto $\overline{D^{\prime}}_P$.}
\end{theorem}

\medskip

In fact, here it is sufficient for the function $K^T_{\mu}(z,z_0)$
to have a dominant of finite mean oscillation in a neighborhood of
every point $z_0\in\partial D$.

\medskip

\begin{corollary}\label{corOSKRSS6.6.2F} {\it In particular, the
conclusion of Theorem \ref{thOSKRSS101F} holds if
\begin{equation}\label{eqOSKRSS6.6.3F}
\overline{\lim\limits_{\varepsilon\to0}}\ \ \
\dashint_{B(z_0,\varepsilon)}K^T_{\mu}(z,z_0)\ dm(z)\ <\
\infty\qquad\forall\ z_0\in\partial D\ .\end{equation}}
\end{corollary}

\begin{theorem}\label{thOSKRSS102F} {\it Let $D$ and $D^{\prime}$ be
bounded finitely connected domains in ${\Bbb C}$ and $f:D\to
D^{\prime}$ be a homeomorphic solution in the class $W^{1,1}_{\rm
loc}$ of the Beltrami equation (\ref{eqBeltrami}) with the condition
\begin{equation}\label{eqOSKRSS10.336aF}\int\limits_{\varepsilon<|z-z_0|<\varepsilon_0}K^T_{\mu}(z,z_0)\
\frac{dm(z)}{|z-z_0|^2}\ =\
o\left(\left[\log\frac{1}{\varepsilon}\right]^2\right)\qquad\forall\
z_0\in\partial D\ .\end{equation} Then $f$ can be extended to a
homeomorphism of $\overline{D}_P$ onto $\overline{D^{\prime}}_P$.}
\end{theorem}

\medskip

\begin{remark}\label{rmOSKRSS200F} Condition
(\ref{eqOSKRSS10.336aF}) can be replaced by the weaker condition
\begin{equation}\label{eqOSKRSS10.336bF}
\int\limits_{\varepsilon<|z-z_0|<\varepsilon_0}\frac{K^T_{\mu}(z,z_0)\
dm(z)}{\left(|z-z_0|\ \log{\frac{1}{|z-z_0|}}\right)^2}\ =\
o\left(\left[\log\log\frac{1}{\varepsilon}\right]^2\right)\qquad\forall\
z_0\in\partial D\ .\end{equation} In general, here we are able to
give a number of other conditions of logarithmic type. In
particular, condition (\ref{eqOSKRSS100dF}), thanking to Theorem
\ref{t:10.1F}, can be replaced by the weaker condition
\begin{equation}\label{eqOSKRSS10.336hF} k^T_{z_0}(r)\ =\ O
\left(\log\frac{1}{r}\log\,\log\frac{1}{r}\right)\ .\end{equation}
\end{remark}

Finally, we complete the series of criteria with the following
integral condition.

\begin{theorem}\label{thOSKRSS103F} {\it Let $D$ and $D^{\prime}$ be
bounded finitely connected domains in ${\Bbb C}$ and $f:D\to
D^{\prime}$ be a homeomorphic solution in the class $W^{1,1}_{\rm
loc}$ of the Beltrami equation (\ref{eqBeltrami}) with the condition
\begin{equation}\label{eqOSKRSS10.36bF}
\int\limits_{D\cap
B(z_0,\varepsilon_0)}\Phi_{z_0}\left(K^T_{\mu}(z,z_0)\right)\ dm(z)\
<\ \infty\qquad\forall\ z_0\in\partial D\end{equation} for
$\varepsilon_0=\varepsilon(z_0)>0$ and a nondecreasing convex
function $\Phi_{z_0}:[0,\infty)\to[0,\infty)$ with
\begin{equation}\label{eqOSKRSS10.37bF}
\int\limits_{\delta_{0}}^{\infty}\frac{d\tau}{\tau\Phi_{z_0}^{-1}(\tau)}\
=\ \infty\end{equation} for $\delta_0=\delta(z_0)>\Phi_{z_0}(0)$.
Then $f$ is extended to a homeomorphism of $\overline{D}_P$ onto
$\overline{D^{\prime}}_P$.}
\end{theorem}

\medskip

\begin{corollary}\label{corOSKRSS6.6.3F} {\it In particular, the
conclusion of Theorem \ref{thOSKRSS103F} holds if
\begin{equation}\label{eqOSKRSS6.6.6F} \int\limits_{D\cap
B(z_0,\varepsilon_0)}e^{\alpha_0 K^T_{\mu}(z,z_0)}\ dm(z)\ <\
\infty\qquad\forall\ z_0\in\partial D\end{equation} for some
$\varepsilon_0=\varepsilon(z_0)>0$ and $\alpha_0=\alpha(z_0)>0$. }
\end{corollary}

\medskip

\begin{remark}\label{rmOSKRSS200000F} Note that condition (\ref{eqOSKRSS10.37bF})
is not only sufficient but also necessary for a continuous extension
to the boundary of all direct mappings $f$ with integral
restrictions of type (\ref{eqOSKRSS10.36bF}), see, e.g., Theorem 5.1
and Remark 5.1 in \cite{KR$_3$}. Recall also that condition
(\ref{eqOSKRSS10.37bF}) is equivalent to each of conditions
(\ref{eq333Y})--(\ref{eq333D}).\end{remark}


\section{Regular solutions for the Dirichlet problem}

Recall that a mapping $f:D\to{\Bbb C}$ is called {\bf discrete} if
the pre-image  $f^{-1}(y)$ of every point $y\in{\Bbb C}$ consists of
isolated points and {\bf open} if the image of every open set
$U\subseteq D$ is open in ${\Bbb C}$.

\medskip

For $\varphi(P)\not\equiv{\rm const}$, $P\in E_D$, a {\bf regular
solution} of Dirichlet problem (\ref{eq1002P}) for Beltrami equation
(\ref{eqBeltrami}) is a continuous discrete open mapping
$f:D\to{\Bbb C}$ of the Sobolev class $W_{\rm loc}^{1,1}$ with its
Jacobian
\begin{equation}\label{eq1003}J_f(z)\ =\ |f_z|^2-|f_{\overline{z}}|^2\ \neq\ 0
\quad\quad{\text{a.e.}}\end{equation} satisfying (\ref{eqBeltrami})
a.e. and condition (\ref{eq1002P}) for all prime ends of the domain
$D$. For $\varphi(P)\equiv c\in{\Bbb R}$, $P\in E_D$, a regular
solution of the problem is any constant function $f(z)=c+ic'$,
$c'\in{\Bbb R}$.


\begin{theorem}\label{thKPRS1001} {\it Let $D$ be a bounded simply
connected domain in $\Bbb C$ and let $\mu:{D}\to {\Bbb D}$ be a
measurable function with $K_{\mu}\in L^{1}_{\mathrm loc}$ and,
moreover,
\begin{equation}\label{eq1008}\int\limits_{0}^{\delta(z_0)}\frac{dr}{||K^T_{\mu}||(z_0,\,r)}\ =\ \infty
\qquad\forall\ z_{0}\in\overline{D}\end{equation} for some
$0<\delta(z_0)<d(z_0)={\sup\limits_{z\in {D}}|z-z_0|}$ and
$$||K^T_\mu||(z_0,\,r)\ :=\ \int\limits_{S(z_0,\,r)}K^T_{\mu}(z, z_0)\ ds\ .$$
Then the Beltrami equation (\ref{eqBeltrami}) has a regular solution
$f$ of the Dirichlet problem (\ref{eq1002P}) for every continuous
function $\varphi:E_D\to{\Bbb R}$.}
\end{theorem}

\medskip

Here and later on, we set that $K^T_{\mu}$ is equal to zero outside
of the domain $D$.

\medskip

\begin{corollary}\label{corKPR1002} {\it Let $D$ be a bounded simply
connected domain in $\Bbb C$ and let $\mu:{D}\to {\Bbb D}$ be a
measurable function such that
\begin{equation}\label{eqOSKRSS100dFD} k^T_{z_{0}}(\varepsilon)\ =\ O{\left(\log\frac{1}{\varepsilon}\right)}
\qquad\mbox{as}\ \ \ \varepsilon\to0 \qquad\forall\
z_{0}\in\overline{D}\end{equation}  where $k^T_{z_{0}}(\varepsilon)$
is the average of the function $K^T_{\mu}(z,z_0)$ over the circle
$S(z_{0},\,\varepsilon)$.

Then the Beltrami equation (\ref{eqBeltrami}) has a regular solution
$f$ of the Dirichlet problem (\ref{eq1002P}) for every continuous
function $\varphi:E_D\to{\Bbb R}$.}
\end{corollary}

\begin{remark}\label{rem2} In particular, the conclusion of Theorem \ref{thKPRS1001} holds if
\begin{equation}\label{eqDIR6**} K^T_{\mu}(z,z_0)=O\left(\log\frac{1}{|z-z_0|}\right)\qquad{\rm
as}\quad z\to z_0\quad\forall\ z_0\in\overline{D}\
.\end{equation}\end{remark}

\medskip

{\it Proof of Theorem \ref{thKPRS1001}.} First of all note that
$E_D$ cannot consist of a single prime end. Indeed, all rays going
from a point $z_0\in D$ to $\infty$ intersect $\partial D$ because
the domain $D$ is bounded, see, e.g., Proposition 2.3 in \cite{RSal}
or Proposition 13.3 in \cite{MRSY}. Thus, $\partial D$ contains more
than one point and by the Riemann theorem, see, e.g., II.2.1 in
\cite{Goluzin}, $D$ can be mapped onto the unit disk ${\Bbb D}$ with
a conformal mapping $R$. However, then by the Caratheodory theorem
there is one-to-one correspondence between elements of $E_D$ and
points of the unit circle $\partial {\Bbb D}$, see, e.g., Theorem
9.6 in \cite{CL}.

Let $F$ be a regular homeomorphic solution of equation
(\ref{eqBeltrami}) in the class $W_{\rm loc}^{1,1}$ which exists in
view of condition (\ref{eq1008}), see, e.g., Theorem 5.4 in paper
\cite{RSY$_3$} or Theorem 11.10 in monograph \cite{MRSY}.

Note that the domain $D^*=F(D)$ is simply connected in
$\overline{\Bbb C}$, see, e.g., Lemma 5.3 in \cite{IR} or Lemma 6.5
in \cite{MRSY}. Let us assume that $\partial D^*$ in $\overline{\Bbb
C}$ consists of the single point $\{ \infty\}$. Then $\overline{\Bbb
C}\setminus D^*$ also consists of the single point $\infty$, i.e.,
$D^*=\Bbb C$, since if there is a point $\zeta_0\in\Bbb C$ in
$\overline{\Bbb C}\setminus D^*$, then, joining it and any point
$\zeta_*\in D^*$ with a segment of a straight line, we find one more
point of $\partial D^*$ in $\Bbb C$, see, e.g., again Proposition
 2.3 in \cite{RSal} or Proposition 13.3 in \cite{MRSY}. Now, let
$\Bbb D^*$ denote the exterior of the unit disk $\Bbb D$ in $\Bbb C$
and let $\kappa(\zeta)=1/\zeta$, $\kappa(0)=\infty$,
$\kappa(\infty)=0$. Consider the mapping $F_*=\kappa\circ
F:\widetilde D\to \Bbb D_0$ where $\widetilde D=F^{-1}(\Bbb D^*)$
and $\Bbb D_0=\Bbb D\setminus \{ 0\}$ is the punctured unit disk. It
is clear that $F_*$ is also a regular homeomorphic solution of
Beltrami equation (\ref{eqBeltrami}) in the class $W_{\rm
loc}^{1,1}$ in the bounded two--connected domain $\widetilde D$
because the mapping $\kappa$ is conformal. By Theorem \ref{t:10.1F}
there is a one--to--one correspondence between elements of $E_D$ and
$0$. However, it was shown above that $E_D$ cannot consists of a
single prime end. This contradiction disproves the above assumption
that $\partial D^*$ consists of a single point in $\overline{\Bbb
C}$.

Thus, by the Riemann theorem $D^*$ can be mapped onto the unit disk
${\Bbb D}$ with a conformal mapping $R_*$. Note that the function
$g:=R_*\circ F$ is again a regular homeomorphic solution in the
Sobolev class $W_{\rm loc}^{1,1}$ of Beltrami equation
(\ref{eqBeltrami}) which maps $D$ onto $\Bbb D$. By Theorem
\ref{t:10.1F} the mapping $g$ admits an extension to a
ho\-meo\-mor\-phism $g_*:{\overline D}_P\to\overline{\Bbb D}$.

We find a regular solution of the initial Dirichlet problem
(\ref{eq1002P}) in the form $f=h\circ g$ where $h$ is a holomorphic
function in $\Bbb D$ with the boundary condition
$$\lim\limits_{z\to\zeta}\,{\rm Re}\,h(z)\ =\
\varphi(g_*^{-1}(\zeta))\qquad \forall\ \zeta\in\partial{\Bbb D}\
.$$ Note that we have from the right hand side a continuous function
of the variable $\zeta$.

As known, the analytic function $h$ can be reconstructed in ${\Bbb
D}$ through its real part on the boundary up to a pure imaginary
additive constant with the Schwartz formula, see, e.g., \S\ 8,
Chapter III, Part 3 in \cite{HuCo},
$$h(z)\ =\ \frac{1}{2\pi i}\int\limits_{|\zeta|=1}\varphi\circ
g_*^{-1}(\zeta)\cdot\frac{\zeta+z}{\zeta-z}\cdot\frac{d\zeta}{\zeta}\
.$$ It is easy to see that the function $f=h\circ g$ is a desired
regular solution of the Dirichlet problem (\ref{eq1002P}) for
Beltrami equation (\ref{eqBeltrami}). $\ \Box$

\medskip

Applying Lemma 2.2 in \cite{RS}, see also Lemma 7.4 in \cite{MRSY},
we obtain the following general lemma immediately from Theorem
\ref{thKPRS1001}.

\medskip

\begin{lemma}\label{lemKPRS1000D} {\it Let $D$ be a bounded simply
connected domain in $\Bbb C$ and $\mu:{D}\to {\Bbb D}$ be a
measurable function with $K_{\mu}\in L^{1}({D})$. Suppose that, for
every $z_0\in\overline{D}$ and every small enough
$\varepsilon_0<d(z_0):=\sup\limits_{z\in D}|z-z_0|$, there is a
family of measurable functions $\psi_{z_0,\,\varepsilon,
\varepsilon_0}:(0,\infty)\to[0,\infty]$,
$\varepsilon\in(0,\,\varepsilon_0)$ such that
\begin{equation}\label{eqKPRS1000}
0\ <\ I_{z_0,\varepsilon_0}(\varepsilon)\ :=\
\int\limits_{\varepsilon}^{\varepsilon_0}
\psi_{z_0,\,\varepsilon,\varepsilon_0}(t)\ dt\ < \
\infty\qquad\forall\ \varepsilon\in(0,\,\varepsilon_0)\end{equation}
and
\begin{equation}\label{eqKPRS1000a}\int\limits_{D(z_0,\,\varepsilon,\,\varepsilon_0)}
K^T_{\mu}(z,z_0)\cdot\psi^{2}_{z_0,\,\varepsilon,
\varepsilon_0}\left(|z-z_0|\right)\, dm(z)\ =\ o(I_{z_0,
\varepsilon_0}^{2}(\varepsilon))\qquad \mbox{as}\ \ \
\varepsilon\to0\end{equation} where
$D(z_0,\,\varepsilon,\,\varepsilon_0)=\{z\in
D:\varepsilon<|z-z_0|<\varepsilon_0\}$. Then the Beltrami equation
(\ref{eqBeltrami}) has a regular solution $f$ of the Dirichlet
problem (\ref{eq1002P}) for every continuous function
$\varphi:E_D\to{\Bbb R}$.} \end{lemma}


\begin{remark}\label{rmKR2.9D}
In fact, it is sufficient here to request instead of the condition
$K_{\mu}\in L^1(D)$ only a local integrability of $K_{\mu}$ in the
domain $D$ and the condition $||K_{\mu}||(z_0,r)\ne\infty$ for a.e.
$r\in (0,\varepsilon_0)$ at all $z_0\in\partial D$.
\end{remark}


By Lemma \ref{lemKPRS1000D} with the choice
$\psi_{z_0,\,\varepsilon}(t)\equiv 1/\left(t\log\frac{1}{t}\right)$
we obtain the following result, see also Lemma \ref{lem13.4.2}.


\begin{theorem}\label{thKPRS1000} {\it Let $D$ be a bounded simply
connected domain in $\Bbb C$ and let $\mu:{D}\to {\Bbb D}$ be a
measurable function with $K_{\mu}\in L^{1}_{\mathrm loc}$ and
\begin{equation}
\label{eq1007}{K^T_{\mu}(z,z_0)\leqslant Q_{z_0}(z)\in{\rm FMO}(z_0)}\qquad\forall\ z_0\in\overline{D} \ .
\end{equation}
Then the Beltrami equation (\ref{eqBeltrami}) has a regular solution
$f$ of the Dirichlet problem (\ref{eq1002P}) for every continuous
function $\varphi:E_D\to{\Bbb R}$.} \end{theorem}

\begin{remark}\label{rm555} In particular, the hypotheses and the conclusion of Theorem
\ref{thKPRS1000} hold if either $Q_{z_0}\in{\rm BMO}_{\rm loc}$ or
$Q_{z_0}\in{\rm W}^{1,2}_{\rm loc}$ because $W^{\,1,2}_{\rm loc}
\subset {\rm VMO}_{\rm loc}$, see, e.g., \cite{BN}.
\end{remark}



By Corollary \ref{FMO_cor2.1} we obtain from Theorem
\ref{thKPRS1000} the following statement.


\begin{corollary}\label{corKPR1001} {\it Let $D$ be a bounded simply
connected domain in $\Bbb C$ and let $\mu:{D}\to {\Bbb D}$ be a
measurable function  with $K_{\mu}\in L^{1}_{\mathrm loc}$ such that
\begin{equation}\label{eq1007b}\limsup\limits_{\varepsilon\to 0}\
\dashint_{B(z_0,\,\varepsilon)}K^T_{\mu}(z,z_0)\ dm(z)\ <\
\infty\qquad\forall\ z_{0}\in{\overline{D}}\ .\end{equation} Then
the Beltrami equation (\ref{eqBeltrami}) has a regular solution $f$
of the Dirichlet problem (\ref{eq1002P}) for every continuous
function $\varphi:E_D\to{\Bbb R}$. }
\end{corollary}

\begin{remark}\label{corKPR1000} {\it In particular, by (\ref{eqConnect}) the
conclusion of Theorem \ref{thKPRS1000} holds if
\begin{equation}\label{eq1007a}
K_{\mu}(z)\ \leqslant\ Q(z)\ \in\ {\rm BMO}({\overline{D}})\ .
\end{equation}}
\end{remark}


The next statement follows from Lemma \ref{lemKPRS1000D} under the
choice $\psi(t)=1/t$, see also Remark \ref{rmKR2.9D}.


\begin{theorem}\label{thOSKRSS102FD} {\it Let $D$ be a bounded simply
connected domain in $\Bbb C$ and let $\mu:{D}\to {\Bbb D}$ be a
measurable function such that
\begin{equation}\label{eqOSKRSS10.336aFD}\int\limits_{\varepsilon<|z-z_0|<\varepsilon_0}K^T_{\mu}(z,z_0)\
\frac{dm(z)}{|z-z_0|^2}\ =\
o\left(\left[\log\frac{1}{\varepsilon}\right]^2\right)\qquad\forall\
z_0\in\overline D\ .\end{equation} Then the Beltrami equation
(\ref{eqBeltrami}) has a regular solution $f$ of the Dirichlet
problem (\ref{eq1002P}) for every continuous function
$\varphi:E_D\to{\Bbb R}$.} \end{theorem}


\begin{remark}\label{rmOSKRSS200FD} Similarly, choosing in Lemma
 \ref{lemKPRS1000D} $\psi(t)=1/(t\log
1/t)$ instead of $\psi(t)=1/t$, we obtain that condition
(\ref{eqOSKRSS10.336aFD}) can be replaced by the condition
\begin{equation}\label{eqOSKRSS10.336bFD}
\int\limits_{\varepsilon<|z-z_0|<\varepsilon_0}\frac{K^T_{\mu}(z,
z_0)\,dm(z)}{\left(|z-z_0|\ \log{\frac{1}{|z-z_0|}}\right)^2}\ =\
o\left(\left[\log\log\frac{1}{\varepsilon}\right]^2\right)\qquad\forall\
z_0\in\overline D\ .\end{equation} Here we are able to give a number
of other conditions of a logarithmic type. In particular, condition
(\ref{eqOSKRSS100dFD}), thanking to Theorem \ref{thKPRS1001}, can be
replaced by the weaker condition
\begin{equation}\label{eqOSKRSS10.336hFD} k^T_{z_0}(r)=O
\left(\log\frac{1}{r}\log\,\log\frac{1}{r}\right).\end{equation}
\end{remark}

Finally, by Theorem \ref{thKPRS1001}, applying also Theorem 3.1 in
\cite{RSY}, we come to the following result.


\begin{theorem}\label{thKPRS1002} {\it Let $D$ be a bounded simply
connected domain in $\Bbb C$ and let $\mu:{D}\to {\Bbb D}$ be a
measurable function with $K_{\mu}\in L^{1}_{loc}$ and
\begin{equation}\label{eq1009}
\int\limits_{D\cap
B(z_0,\varepsilon_0)}\Phi_{z_0}(K^T_{\mu}(z,z_0))\ dm(z)\ <\ \infty
\qquad\forall\ z_0\in\overline D
\end{equation} for $\varepsilon_0=\varepsilon(z_0)>0$ and a nondecreasing convex function $\Phi_{z_0}:[0,\infty)\to[0,\infty)$
with
\begin{equation}
\label{eq1010D}\int\limits_{\delta_0}^{\infty}\frac{d\tau}{\tau\Phi_{z_0}^{-1}(\tau)}\
=\ \infty
\end{equation}
for $\delta_0=\delta(z_0)>\Phi_{z_0}(0)$. Then the Beltrami equation
(\ref{eqBeltrami}) has a regular solution $f$ of the Dirichlet
problem (\ref{eq1002P}) for every continuous function
$\varphi:E_D\to{\Bbb R}$.}
\end{theorem}


\begin{remark}\label{rmkKPRS1000} Recall that condition (\ref{eq1010D})
is equivalent to each of conditions (\ref{eq333Y})--(\ref{eq333D}).
Moreover, condition (\ref{eq1010D}) is not only sufficient but also
ne\-ces\-sa\-ry to have a regular solution of the Dirichlet problem
(\ref{eq1002P}) for every Beltrami equation (\ref{eqBeltrami}) with
integral restriction (\ref{eq1009}) for every continuous function
$\varphi:E_D\to{\Bbb R}$. Indeed, by the Stoilow theorem on
representation of discrete open mappings, see, e.g., \cite{Sto},
every regular solution $f$ of the Dirichlet problem (\ref{eq1002P})
for Beltrami equation (\ref{eqBeltrami}) with $K_{\mu}\in
L^{1}_{loc}$ can be represented in the form of composition
$f=h\circ{F}$ where $h$ is a holomorphic function and $F$ is a
regular homeomorphic solution of (\ref{eqBeltrami}) in the class
$W_{\rm loc}^{1,1}$. Thus, by Theorem 5.1 in \cite{RSY13} on the
nonexistence of regular homeomorphic solutions of (\ref{eqBeltrami})
in the class $W_{\rm loc}^{1,1}$, if (\ref{eq1010D}) fails, then
there is a measurable function $\mu:{D}\to {\Bbb D}$ satisfying
integral condition (\ref{eq1009}) for which Beltrami equation
(\ref{eqBeltrami}) has no regular solution of the Dirichlet problem
(\ref{eq1002P}) for any nonconstant continuous function
$\varphi:E_D\to{\Bbb R}$.
\end{remark}


\begin{corollary}\label{corKPR1003} {\it Let $D$ be a bounded simply
connected domain in $\Bbb C$ and let $\mu:{D}\to {\Bbb D}$ be a
measurable function with $K_{\mu}\in L^{1}_{loc}$ and
\begin{equation}\label{eq1011}
\int\limits_{D\cap B(z_0,\varepsilon_0)}e^{\alpha_0
K^T_{\mu}(z,z_0)}\ dm(z)\ <\ \infty\qquad\forall\ z_0\in\overline{D}
\end{equation}  for some
$\varepsilon_0=\varepsilon(z_0)>0$ and $\alpha_0=\alpha(z_0)>0$.
Then the Beltrami equation (\ref{eqBeltrami}) has a regular solution
$f$ of the Dirichlet problem (\ref{eq1002P}) for every continuous
function $\varphi:E_D\to{\Bbb R}$.}
\end{corollary}

\bigskip

\section{Pseudoregular solutions in multiply connected domains}

As it was first noted by Bogdan Bojarski, see, e.g., \S\ 6 of
Chapter 4 in \cite{Vekua}, the Dirichlet problem for the Beltrami
equations, generally speaking, has no regular solution in the class
of continuous (single--valued) in ${\Bbb C}$ functions with
generalized derivatives in the case of multiply connected domains
$D$. Hence the natural question arose: whether solutions exist in
wider classes of functions for this case ? It is turned out to be
solutions for this problem can be found in the class of functions
admitting a certain number (related with connectedness of $D$) of
poles at prescribed points. Later on, this number will take into
account the multiplicity of these poles from the Stoilow
representation.


Namely, a {\bf pseudoregular solution} of such a problem is a
continuous (in $\overline{\Bbb C}$) discrete open mapping
$f:D\to\overline{\Bbb C}$ of the Sobolev class $W_{\rm loc}^{1,1}$
(outside of poles) with its Jacobian
$J_f(z)=\left|f_z\right|^2-\left|f_{\overline{z}}\right|^2\ne0$ a.e.
satisfying (\ref{eqBeltrami}) a.e. and the boundary condition
(\ref{eq1002P}).

\medskip

Arguing similarly to the case of simply connected domains and
applying Theorem V.6.2 in \cite{Goluzin} on conformal mappings of
finitely connected domains onto circular domains and also Theorems
4.13 and 4.14 in \cite{Vekua}, we obtain the following result.


\begin{theorem}\label{thKPRS1001M} {\it Let $D$ be a bounded $m-$connected domain in $\Bbb
C$ with nondegenerate boundary components, $k\geqslant m-1$ and
$\mu:{D}\to {\Bbb D}$ be a measurable function with $K_{\mu}\in
L^{1}_{\mathrm loc}$ and
\begin{equation}\label{eq1008M}\int\limits_{0}^{\delta(z_0)}\frac{dr}{||K^T_{\mu}||(z_0,\,r)}\ =\ \infty
\qquad\forall\ z_{0}\in\overline{D}\end{equation} for some
$0<\delta(z_0)<d(z_0)={\sup\limits_{z\in {D}}|z-z_0|}$ and
$$||K^T_\mu||(z_0,\,r)\ :=\ \int\limits_{S(z_0,\,r)}K^T_{\mu}(z, z_0)\ ds\ .$$
Then the Beltrami equation (\ref{eqBeltrami}) has a pseudoregular
solution $f$ of the Dirichlet problem (\ref{eq1002P}) with $k$ poles
at prescribed points in $D$ for every continuous function
$\varphi:E_D\to{\Bbb R}$.}
\end{theorem}


Here, as before, we set $K^T_{\mu}$ to be extended by zero outside
of the domain $D$.


\begin{corollary}\label{corKPR1002M} {\it Let $D$ be a bounded $m-$connected domain in $\Bbb
C$ with nondegenerate boundary components, $k\geqslant m-1$ and
$\mu:{D}\to {\Bbb D}$ be a measurable function with $K_{\mu}\in
L^{1}_{\mathrm loc}$ and
\begin{equation}\label{eqOSKRSS100dFDM} k^T_{z_{0}}(\varepsilon)=O{\left(\log\frac{1}{\varepsilon}\right)}
\qquad\mbox{as}\ \ \ \varepsilon\to0 \qquad\forall\
z_{0}\in\overline{D}\end{equation} where $k^T_{z_{0}}(\varepsilon)$
is the average of the function $K^T_{\mu}(z,z_0)$ over the circle
$S(z_{0},\,\varepsilon)$.

Then the Beltrami equation (\ref{eqBeltrami}) has a pseudoregular
solution $f$ of the Dirichlet problem (\ref{eq1002P}) with $k$ poles
at prescribed points in $D$ for every continuous function
$\varphi:E_D\to{\Bbb R}$.}
\end{corollary}

\begin{remark}\label{rem2p} In particular, the conclusion of Theorem \ref{thKPRS1001M} holds if
\begin{equation}\label{eqDIR6**p} K^T_{\mu}(z,z_0)=O\left(\log\frac{1}{|z-z_0|}\right)\qquad{\rm
as}\quad z\to z_0\quad\forall\ z_0\in\overline{D}\
.\end{equation}\end{remark}


{\it Proof of Theorem \ref{thKPRS1001M}.} Let $F$ be a regular
solution of equation (\ref{eqBeltrami}) in the class $W_{\rm
loc}^{1,1}$ that exists by condition (\ref{eq1008M}), see, e.g.,
Theorem 5.4 in the paper \cite{RSY$_3$} or Theorem 11.10 in the
monograph \cite{MRSY}. Note that the domain $D^*=F(D)$ is
$m-$connected in $\overline{\Bbb C}$ and there is a natural
one--to--one correspondence between components $\gamma_j$ of $\gamma
=\partial D$ and components $\Gamma_j$ of $\Gamma =\partial D^*$,
$\Gamma_j=C(\gamma_j, F)$ and $\gamma_j=C(\Gamma_j, F^{-1})$,
$j=1,\ldots , m$, see, e.g., Lemma 5.3 in \cite{IR} or Lemma 6.5 in
\cite{MRSY}. Moreover, by Remark \ref{METRIC} every subspace $E_j$
of $E_D$ associated with $\gamma_j$ consists of more than one prime
end, even it is homeomorphic to the unit circle.

Next, no one of $\Gamma_j$, $j=1,\ldots , m$, is degenerated to a
single point. Indeed, let us assume that $\Gamma_{j_0}=\{ \zeta_0\}$
first for some $\zeta_0\in\Bbb C$. Let $r_0\in(0,d_0)$ where
$d_0=\inf\limits_{\zeta\in\Gamma\setminus\Gamma_{j_0}}|\zeta
-\zeta_0|$. Then the punctured disk $D_0=\{ \zeta\in\Bbb C :
0<|\zeta -\zeta_0|<r_0\}$ is in the domain $D^*$ and its boundary
does not intersect $\Gamma\setminus\Gamma_{j_0}$. Set $\widetilde
D=F^{-1}(D_0)$. Then by the construction $\widetilde D\subset D$ is
a 2--connected domain, $\overline{\widetilde
D}\cap\gamma\setminus\gamma_{j_0}=\varnothing$,
$C(\gamma_{j_0},\widetilde F)=\{ \zeta_0\}$ and
$C(\zeta_0,\widetilde F^{-1})=\gamma_{j_0}$ where $\widetilde F$ is
a restriction of the mapping $\widetilde F$ to $\widetilde D$.
However, this contradicts Theorem \ref{t:10.1F} because, as it was
noted above, $E_{j_0}$ contains more than one prime end.

Now, let assume that $\Gamma_{j_0}=\{ \infty\}$. Then the component
of $\overline{\Bbb C}\setminus D^*$ associated with $\Gamma_{j_0}$,
see Lemma 5.1 in \cite{IR} or Lemma 6.3 in \cite{MRSY}, is also
consists of the single point $\infty$ because if the interior of
this component is not empty, then choosing there an arbitrary point
$\zeta_0$ and joining it with a point $\zeta_*\in D^*$ by a segment
of a straight line we would find one more point in $\Gamma_{j_0}$,
see, e.g., Proposition 2.3 in \cite{RSal} or Proposition 13.3 in
\cite{MRSY}.

Thus, applying if it is necessary an additional stretching
(conformal mapping), with no loss of generality we may assume that
$D^*$ contains the exteriority $\Bbb D_*$ of the unit disk $\Bbb D$
in $\Bbb C$. Set $\kappa(\zeta)=1/\zeta$, $\kappa(0)=\infty$,
$\kappa(\infty)=0$. Consider the mapping $F_*=\kappa\circ
F:\widetilde D\to \Bbb D_0$ where $\widetilde D=F^{-1}(\Bbb D_*)$
and $\Bbb D_0=\Bbb D\setminus \{ 0\}$ is the punctured unit disk. It
is clear that $F_*$ is also a homeomorphic solution of Beltrami
equation (\ref{eqBeltrami}) of the class $W_{\rm loc}^{1,1}$ in
2--connected domain $\widetilde D$ because the mapping $\kappa$ is
conformal. Consequently, by Theorem
 \ref{t:10.1F} elements of $E_{j_0}$ should be in a one--to--one
 correspondence with $0$. However, it was already noted, $E_{j_0}$
cannot consists of a single prime end. The obtained contradiction
disproves the assumption that $\Gamma_{j_0}=\{ \infty\}$.

Thus, by Theorem V.6.2 â \cite{Goluzin}, see also Remark 1.1 in
\cite{KPR$_*$}, $D^*$ can be mapped with a conformal mapping $R_*$
onto a bounded circular domain ${\Bbb D}^*$ whose  boundary consists
of mutually disjoint circles. Note that the function $g:=R_*\circ F$
is again a regular homeomorphic solution in the Sobolev class
$W_{\rm loc}^{1,1}$ for Beltrami equation (\ref{eqBeltrami}) that
maps $D$ onto $\Bbb D^*$. By Theorem \ref{t:10.1F} the mapping $g$
admits an extension to a homeomorphism $g_*:{\overline
D}_P\to\overline{\Bbb D^*}$.

Let us find a solution of the initial Dirichlet problem
(\ref{eq1002P}) in the form $f=h\circ g$ where $h$ is a meromorphic
function in $\Bbb D^*$ with the boundary condition
\begin{equation}\label{BOUNDARY}\lim\limits_{z\to\zeta}\,{\rm
Re}\,h(z)\ =\ \varphi(g_*^{-1}(\zeta))\qquad \forall\
\zeta\in\partial{\Bbb D^*} \end{equation} and $k\geqslant m-1$ poles
corresponding under the mapping $g$ to those at prescribed points in
$D$. Note that the function from the right hand side in
(\ref{BOUNDARY}) is continuous in the variable $\zeta$. Thus, such a
function $h$ exists by Theorems 4.13 and 4.14 in \cite{Vekua}. It is
clear that the function $f$ associated with $h$ is by the
construction a desired pseudoregular solution of the Dirichlet
problem (\ref{eq1002P}) for Beltrami equation (\ref{eqBeltrami}). $\
\Box$

\medskip

Applying Lemma 2.2 in \cite{RS}, see also Lemma 7.4 in \cite{MRSY},
we obtain immediately from Theorem \ref{thKPRS1001M} the next lemma.


\begin{lemma}\label{lemKPRS1000DM} {\it Let $D$ be a bounded $m-$connected domain in $\Bbb
C$ with nondegenerate boundary components,  $k\geqslant m-1$ and
$\mu:{D}\to {\Bbb D}$ be a measurable function with $K_{\mu}\in
L^{1}({D})$. Suppose that, for every $z_0\in\overline{D}$ and every
small enough $0<\varepsilon_0<d(z_0):=\sup\limits_{z\in D}|z-z_0|$,
there is a family of measurable functions
$\psi_{z_0,\,\varepsilon,\varepsilon_0}:(0,\infty)\to[0,\infty]$,
$\varepsilon\in(0,\,\varepsilon_0)$ such that
\begin{equation}\label{eqKPRS1000M}
0\ <\ I_{z_0,\varepsilon_0}(\varepsilon)\ :=\
\int\limits_{\varepsilon}^{\varepsilon_0}
\psi_{z_0,\,\varepsilon,\varepsilon_0}(t)\ dt\ < \
\infty\qquad\forall\ \varepsilon\in(0,\,\varepsilon_0)\end{equation}
and
\begin{equation}\label{eqKPRS1000aM}\int\limits_{\varepsilon<|z-z_0|<\varepsilon_0}
K^T_{\mu}(z,z_0)\cdot\psi^{2}_{z_0,\,\varepsilon,\varepsilon_0}\left(|z-z_0|\right)\,
dm(z)\ =\ o(I_{z_0,\varepsilon_0}^{2}(\varepsilon))\ \ \ \ \ \ \
\mbox{as}\ \ \ \varepsilon\to0\ .\end{equation} Then the Beltrami
equation (\ref{eqBeltrami}) has a pseudoregular solution $f$ of the
Dirichlet problem (\ref{eq1002P}) with $k$ poles at prescribed
points in $D$ for every continuous function $\varphi:E_D\to{\Bbb
R}$.}
\end{lemma}


\begin{remark}\label{rmKR2.9DM}
In fact, here it is sufficient to assume instead of the condition
$K_{\mu}\in L^1(D)$ the local integrability of $K_{\mu}$ in the
domain $D$ and the condition $||K_{\mu}||(z_0,r)\ne\infty$ for a.e.
$r\in (0,\varepsilon_0)$ and all $z_0\in\partial D$.
\end{remark}


By Lemma \ref{lemKPRS1000DM} with the choice
$\psi_{z_0,\,\varepsilon}(t)\equiv 1/t\log\frac{1}{t}$ we obtain the
following result, see also  Lemma  \ref{lem13.4.2}.


\begin{theorem}\label{thKPRS1000M} {\it Let $D$ be a bounded $m-$connected domain in $\Bbb
C$ with nondegenerate boundary components, $k\geqslant m-1$ and
$\mu:{D}\to {\Bbb D}$ be a measurable function  with $K_{\mu}\in
L^{1}({D})$ such that
\begin{equation}\label{eq1007M} {K^T_{\mu}(z,z_0)\ \leqslant\ Q_{z_0}(z)\ \in\ {\rm FMO}(z_0)}\qquad\forall\ z_0 \in\overline{D} \ .
\end{equation} Then the Beltrami equation (\ref{eqBeltrami}) has a
pseudoregular solution $f$ of the Dirichlet problem (\ref{eq1002P})
with $k$ poles at prescribed points in $D$ for every continuous
function $\varphi:E_D\to{\Bbb R}$.}
\end{theorem}

\begin{remark}\label{rm555p} In particular, the conclusion of Theorem
\ref{thKPRS1000M} holds if either $Q_{z_0}\in{\rm BMO}_{\rm loc}$ or
$Q_{z_0}\in{\rm W}^{1,2}_{\rm loc}$ because $W^{\,1,2}_{\rm loc}
\subset {\rm VMO}_{\rm loc}$, see, e.g., \cite{BN}.
\end{remark}


By Corollary  \ref{FMO_cor2.1} we have the next consequence of
Theorem \ref{thKPRS1000M}:


\begin{corollary}\label{corKPR1001M} {\it Let $D$ be a bounded $m-$connected domain in $\Bbb
C$ with nondegenerate boundary components, $k\geqslant m-1$ and
$\mu:{D}\to {\Bbb D}$ be a measurable function  with $K_{\mu}\in
L^{1}({D})$ such that
\begin{equation}\label{eq1007Mb}
\limsup\limits_{\varepsilon\to 0}\
\dashint_{B(z_0,\,\varepsilon)}K^T_{\mu}(z,z_0)\ dm(z)\ <\
\infty\qquad\forall\ z_{0}\ \in\ {\overline{D}}\ .\end{equation}
Then the Beltrami equation (\ref{eqBeltrami}) has a pseudoregular
solution $f$ of the Dirichlet problem (\ref{eq1002P}) with $k$ poles
at prescribed points in $D$ for every continuous function
$\varphi:E_D\to{\Bbb R}$.}
\end{corollary}


\begin{remark}\label{corKPR1000M} {\it In particular, by (\ref{eqConnect}) the
conclusion of Theorem \ref{thKPRS1000M} holds if
\begin{equation}\label{eq1007Ma} K_{\mu}(z)\ \leqslant\ Q(z)\ \in\ {\rm
BMO}({\overline{D}})\end{equation}.}
\end{remark}


The following statement follows from Lemma \ref{lemKPRS1000DM}
through the choice $\psi(t)=1/t$, see also Remark \ref{rmKR2.9DM}.


\begin{theorem}\label{thOSKRSS102FDM} {\it Let $D$ be a bounded $m-$connected domain in $\Bbb
C$ with nondegenerate boundary components, $k\geqslant m-1$ and
$\mu:{D}\to {\Bbb D}$ be a measurable function such that
\begin{equation}\label{eqOSKRSS10.336aFDM}\int\limits_{\varepsilon<|z-z_0|<\varepsilon_0}K^T_{\mu}(z,z_0)\
\frac{dm(z)}{|z-z_0|^2}\ =\
o\left(\left[\log\frac{1}{\varepsilon}\right]^2\right)\qquad\forall\
z_0\in\overline D\ .\end{equation} Then the Beltrami equation
(\ref{eqBeltrami}) has a pseudoregular solution $f$ of the Dirichlet
problem (\ref{eq1002P}) with $k$ poles at prescribed points in $D$
for every continuous function $\varphi:E_D\to{\Bbb R}$.}
\end{theorem}


\begin{remark}\label{rmOSKRSS200FDM} Similarly, choosing in Lemma
 \ref{lemKPRS1000DM} $\psi(t)=1/(t\log
1/t)$ instead of $\psi(t)=1/t$ we obtain that condition
(\ref{eqOSKRSS10.336aFDM}) can be replaced by the condition
\begin{equation}\label{eqOSKRSS10.336bFDM}
\int\limits_{\varepsilon<|z-z_0|<\varepsilon_0}\frac{K^T_{\mu}(z,z_0)\,dm(z)}{\left(|z-z_0|\
\log{\frac{1}{|z-z_0|}}\right)^2}\ =\
o\left(\left[\log\log\frac{1}{\varepsilon}\right]^2\right)\qquad\forall\
z_0\in\overline D\ .\end{equation} Here we are able to give a number
of other conditions of the logarithmic type. In particular,
condition (\ref{eqOSKRSS100dFDM}), thanking to Theorem
\ref{thKPRS1001M}, can be replaced by the weaker condition
\begin{equation}\label{eqOSKRSS10.336hFDM} k^T_{z_0}(r)=O
\left(\log\frac{1}{r}\log\,\log\frac{1}{r}\right).\end{equation}
\end{remark}

Finally, by Theorem \ref{thKPRS1001M}, applying also Theorem 3.1 in
the paper \cite{RSY}, we come to the following result.


\begin{theorem}\label{thKPRS1002M} {\it Let $D$ be a bounded $m-$connected domain in $\Bbb
C$ with nondegenerate boundary components, $k\geqslant m-1$ and
$\mu:{D}\to {\Bbb D}$ be a measurable function with $K_{\mu}\in
L^{1}_{loc}$ such that
\begin{equation}\label{eq1009M}\int\limits_{D\cap B(z_0,\varepsilon_0)}\Phi_{z_0}(K^T_{\mu}(z,z_0))\ dm(z)\ <\ \infty\end{equation}
for $\varepsilon_0=\varepsilon(z_0)>0$ and a nondecreasing convex
function $\Phi_{z_0}:[0,\infty)\to[0,\infty)$ with
\begin{equation}\label{eq1010M}\int\limits_{\delta_0}^{\infty}\frac{d\tau}{\tau\Phi_{z_0}^{-1}(\tau)}\ =\ \infty\end{equation}
for $\delta_0=\delta(z_0)>\Phi_{z_0}(0)$. Then the Beltrami equation
(\ref{eqBeltrami}) has a pseudoregular solution $f$ of the Dirichlet
problem (\ref{eq1002P}) with $k$ poles at prescribed points in $D$
for every continuous function $\varphi:E_D\to{\Bbb R}$.}
\end{theorem}


Recall that condition (\ref{eq1010M}) is equivalent to every of
conditions (\ref{eq333Y})--(\ref{eq333D}).


\begin{corollary}\label{corKPR1003M} {\it Let $D$ be a bounded $m-$connected domain in $\Bbb
C$ with nondegenerate boundary components, $k\geqslant m-1$ and
$\mu:{D}\to {\Bbb D}$ be a measurable function with $K_{\mu}\in
L^{1}_{loc}$ such that
\begin{equation}\label{eq1011M} \int\limits_{D\cap B(z_0,\varepsilon_0)}e^{\alpha_0
K^T_{\mu}(z,z_0)}\ dm(z)\ <\ \infty\qquad\forall\ z_0\in\overline{D}
\end{equation}  for some
$\varepsilon_0=\varepsilon(z_0)>0$ and $\alpha_0=\alpha(z_0)>0$.

Then the Beltrami equation (\ref{eqBeltrami}) has a pseudoregular
solution $f$ of the Dirichlet problem (\ref{eq1002P}) with $k$ poles
at prescribed points in $D$ for every continuous function
$\varphi:E_D\to{\Bbb R}$.}
\end{corollary}

\bigskip

\section{Multivalent solutions in finitely connected domains}

In  finitely connected domains $D$ in $\Bbb{C}$, in addition to pseudoregular solutions, the Dirichlet problem
(\ref{eqGrUsl}) for the Beltrami equation (\ref{eqBeltrami}) admits multi-valued solutions in the spirit of the
theory of multi-valued analytic functions. We say that a discrete open mapping $f:B(z_0,\varepsilon_0)\to{\Bbb
C}$, where $B(z_0,\varepsilon_0)\subseteq D$, is a {\bf local regular solution of the equation}
(\ref{eqBeltrami}) if $f\in W_{\rm loc}^{1,1}$, $J_f(z)\neq0$ and $f$ satisfies (\ref{eqBeltrami}) a.e. in
$B(z_0,\varepsilon_0)$.

\medskip

The local regular solutions $f:B(z_0,\varepsilon_0)\to{\Bbb C}$ and
$f_*:B(z_*,\varepsilon_*)\to{\Bbb C}$ of the equation
(\ref{eqBeltrami}) will be called extension of each to other if
there is a finite chain of such solutions
$f_i:B(z_i,\varepsilon_i)\to\Bbb{C}$, $i=1,\ldots,m$, that
$f_1=f_0$, $f_m=f_*$ and $f_i(z)\equiv f_{i+1}(z)$ for $z\in
E_i:=B(z_i,\varepsilon_i)\cap
B(z_{i+1},\varepsilon_{i+1})\neq\emptyset$, $i=1,\ldots,m-1$. A
collection of local regular solutions
$f_j:B(z_j,\varepsilon_j)\to{\Bbb C}$, $j\in J$, will be called a
{\bf multivalent solution} of the equation (\ref{eqBeltrami}) in $D$
if the disks $B(z_j,\varepsilon_j)$ cover the whole domain $D$ and
$f_j$ are extensions of each to other through the collection and the
collection is maximal by inclusion. A multi-valued solution of the
equation (\ref{eqBeltrami}) will be called a {\bf multivalent
solution of the Dirichlet problem} (\ref{eqGrUsl}) if $u(z)={\rm
Re}\,f(z)={\rm Re}\,f_{j}(z)$, $z\in B(z_j,\varepsilon_j)$, $j\in
J$, is a single-valued function in $D$ satisfying the condition
$\lim\limits_{z\in\zeta}u(z)=\varphi(\zeta)$ for all $\zeta\to
\partial D$.

\medskip

As it was before, we assume later on that $K^T_{\mu}(\cdot,z_0)$ is
extended by zero outside of the domain $D$.

\medskip

The proof of the existence of multivalent solutions of Dirichlet
problem (\ref{eq1002P}) for Beltrami equation (\ref{eqBeltrami}) in
finitely connected domains is reduced to the Dirichlet problem for
harmonic functions in circular domains, see, e.g., \S\ 3 of Chapter
VI in \cite{Goluzin}.

\begin{theorem}\label{thKPRS1001MM} {\it Let $D$ be a bounded finitely connected domain in $\Bbb
C$ with nondegenerate boundary components and $\mu:{D}\to {\Bbb D}$
be a measurable function with $K_{\mu}\in L^{1}_{\mathrm loc}$ and
\begin{equation}\label{eq1008MM}\int\limits_{0}^{\delta(z_0)}\frac{dr}{||K^T_{\mu}||(z_0,\,r)}\ =\ \infty
\qquad\forall\ z_{0}\in\overline{D}\end{equation} for some
$0<\delta(z_0)<d(z_0)={\sup\limits_{z\in {D}}|z-z_0|}$ and
$$||K^T_\mu||(z_0,\,r)\ :=\ \int\limits_{S(z_0,\,r)}K^T_{\mu}(z, z_0)\ ds\ .$$
Then Beltrami equation (\ref{eqBeltrami}) has a multivalent solution
of Dirichlet problem (\ref{eq1002P}) for every continuous function
$\varphi:E_{D}\to{\Bbb R}$.}
\end{theorem}


{\it Proof of Theorem \ref{thKPRS1001MM}.} Similarly to the first
part of Theorem \ref{thKPRS1001MM}, it is proved that there is a
regular homeomorphic solution $g$ of Beltrami equation
(\ref{eqBeltrami}) mapping the domain $D$ onto a circular domain
$\Bbb D^*$ whose boundary consists of mutually disjoint circles. By
Theorem \ref{t:10.1F} the mapping $g$ admits an extension to a
homeomorphism $g_*:{\overline D}_P\to\overline{\Bbb D^*}$.

As known, in the circular domain $\Bbb D^*$, there is a solution of
the Dirichlet problem
\begin{equation}\label{BOUNDARYMM}\lim\limits_{z\to\zeta}\ u(z)\ =\ \varphi(g_*^{-1}(\zeta))\qquad \forall\
\zeta\in\partial{\Bbb D^*} \end{equation} for harmonic functions
$u$, see, e.g., \S\ 3 of Chapter VI in \cite{Goluzin}. Let
$B_0=B(z_0,r_0)$ is a disk in the domain $D$. Then ${\frak B}_0 =
g(B_0)$ is a simply connected subdomain of the circular domain $\Bbb
D^*$ where there is a conjugate function $v$ determined up to an
additive constant such that $h=u+iv$ is a single--valued analytic
function. The function $h$ can be extended to, generally speaking
multivalent, analytic function $H$ along any path in $\Bbb D^*$
because $u$ is given in the whole domain  $\Bbb D^*$.

Thus, $f=H\circ g$ is a desired multivalent solution of the
Dirichlet problem (\ref{eq1002P}) for Beltrami equation
(\ref{eqBeltrami}).  $\ \Box$

\medskip

The hypotheses of the rest theorems and corollaries below imply the
hypotheses of Theorem \ref{thKPRS1001MM} as it was shown in the
previous section.

\begin{corollary}\label{corKPR1002MM} {\it Let $D$ be a bounded finitely connected domain in $\Bbb
C$ with nondegenerate boundary components and $\mu:{D}\to {\Bbb D}$
be a measurable function with $K_{\mu}\in L^{1}_{\mathrm loc}$ and
\begin{equation}\label{eqOSKRSS100dFDMM} k^T_{z_{0}}(\varepsilon)=O{\left(\log\frac{1}{\varepsilon}\right)}
\qquad\mbox{as}\ \ \ \varepsilon\to0 \qquad\forall\
z_{0}\in\overline{D}\end{equation} where $k^T_{z_{0}}(\varepsilon)$
is the average of the function $K^T_{\mu}(z,z_0)$ over the circle
$S(z_{0},\,\varepsilon)$.

Then Beltrami equation (\ref{eqBeltrami}) has a multivalent solution
of Dirichlet problem (\ref{eq1002P}) for every continuous function
$\varphi:E_{D}\to{\Bbb R}$. }
\end{corollary}

\begin{remark}\label{rem2pM} In particular, the conclusion of Theorem \ref{thKPRS1001MM} holds if
\begin{equation}\label{eqDIR6**pM} K^T_{\mu}(z,z_0)=O\left(\log\frac{1}{|z-z_0|}\right)\qquad{\rm
as}\quad z\to z_0\quad\forall\ z_0\in\overline{D}\
.\end{equation}\end{remark}

\medskip

Applying Lemma 2.2 in \cite{RS}, see also Lemma 7.4 in \cite{MRSY},
we obtain immediately from Theorem \ref{thKPRS1001MM} the next
lemma.


\begin{lemma}\label{lemKPRS1000DMM} {\it Let $D$ be a bounded finitely connected domain in $\Bbb
C$ with nondegenerate boundary components and $\mu:{D}\to {\Bbb D}$
be a measurable function with $K_{\mu}\in L^{1}({D})$. Suppose that,
for every $z_0\in\overline{D}$ and every small enough
$0<\varepsilon_0<d(z_0):=\sup\limits_{z\in D}|z-z_0|$, there is a
family of measurable functions
$\psi_{z_0,\,\varepsilon,\varepsilon_0}:(0,\infty)\to[0,\infty]$,
$\varepsilon\in(0,\,\varepsilon_0)$ such that
\begin{equation}\label{eqKPRS1000MM}
0\ <\ I_{z_0,\varepsilon_0}(\varepsilon)\ :=\
\int\limits_{\varepsilon}^{\varepsilon_0}
\psi_{z_0,\,\varepsilon,\varepsilon_0}(t)\ dt\ < \
\infty\qquad\forall\ \varepsilon\in(0,\,\varepsilon_0)\end{equation}
and
\begin{equation}\label{eqKPRS1000aMM}\int\limits_{\varepsilon<|z-z_0|<\varepsilon_0}
K^T_{\mu}(z,z_0)\cdot\psi^{2}_{z_0,\,\varepsilon,\varepsilon_0}\left(|z-z_0|\right)\,
dm(z)\ =\ o(I_{z_0,\varepsilon_0}^{2}(\varepsilon))\ \ \ \ \ \ \
\mbox{as}\ \ \ \varepsilon\to0\ .\end{equation} Then Beltrami
equation (\ref{eqBeltrami}) has a multivalent solution of Dirichlet
problem (\ref{eq1002P}) for every continuous function
$\varphi:E_{D}\to{\Bbb R}$.}
\end{lemma}


\begin{remark}\label{rmKR2.9DMM}
In fact, here it is sufficient to assume instead of the condition
$K_{\mu}\in L^1(D)$ the local integrability of $K_{\mu}$ in the
domain $D$ and the condition $||K_{\mu}||(z_0,r)\ne\infty$ for a.e.
$r\in (0,\varepsilon_0)$ and all $z_0\in\partial D$.
\end{remark}


By Lemma \ref{lemKPRS1000DMM} with the choice
$\psi_{z_0,\,\varepsilon}(t)\equiv 1/t\log\frac{1}{t}$ we obtain the
following result, see also  Lemma  \ref{lem13.4.2}.


\begin{theorem}\label{thKPRS1000MM} {\it Let $D$ be a bounded finitely connected domain in $\Bbb
C$ with nondegenerate boundary components and $\mu:{D}\to {\Bbb D}$
be a measurable function  with $K_{\mu}\in L^{1}({D})$ such that
\begin{equation}\label{eq1007MM} {K^T_{\mu}(z,z_0)\ \leqslant\ Q_{z_0}(z)\ \in\ {\rm FMO}(z_0)}\qquad\forall\ z_0 \in\overline{D} \ .
\end{equation} Then Beltrami equation (\ref{eqBeltrami}) has a multivalent solution
of Dirichlet problem (\ref{eq1002P}) for every continuous function
$\varphi:E_{D}\to{\Bbb R}$.}
\end{theorem}

\begin{remark}\label{rm555p} In particular, the conclusion of Theorem
\ref{thKPRS1000MM} holds if either $Q_{z_0}\in{\rm BMO}_{\rm loc}$
or $Q_{z_0}\in{\rm W}^{1,2}_{\rm loc}$ because $W^{\,1,2}_{\rm loc}
\subset {\rm VMO}_{\rm loc}$, see, e.g., \cite{BN}.
\end{remark}


By Corollary  \ref{FMO_cor2.1} we have the next consequence of
Theorem \ref{thKPRS1000MM}:


\begin{corollary}\label{corKPR1001MM} {\it Let $D$ be a bounded finitely connected domain in $\Bbb
C$ with nondegenerate boundary components and $\mu:{D}\to {\Bbb D}$
be a measurable function  with $K_{\mu}\in L^{1}({D})$ such that
\begin{equation}\label{eq1007Mb}
\limsup\limits_{\varepsilon\to 0}\
\dashint_{B(z_0,\,\varepsilon)}K^T_{\mu}(z,z_0)\ dm(z)\ <\
\infty\qquad\forall\ z_{0}\ \in\ {\overline{D}}\ .\end{equation}
Then Beltrami equation (\ref{eqBeltrami}) has a multivalent solution
of Dirichlet problem (\ref{eq1002P}) for every continuous function
$\varphi:E_{D}\to{\Bbb R}$. }
\end{corollary}


\begin{remark}\label{corKPR1000MM} {\it In particular, by (\ref{eqConnect}) the
conclusion of Theorem \ref{thKPRS1000MM} holds if
\begin{equation}\label{eq1007Ma} K_{\mu}(z)\ \leqslant\ Q(z)\ \in\ {\rm
BMO}({\overline{D}})\end{equation}.}
\end{remark}


The following statement follows from Lemma \ref{lemKPRS1000DMM}
through the choice $\psi(t)=1/t$, see also Remark \ref{rmKR2.9DMM}.


\begin{theorem}\label{thOSKRSS102FDMM} {\it Let $D$ be a bounded finitely connected domain in $\Bbb
C$ with nondegenerate boundary components and $\mu:{D}\to {\Bbb D}$
be a measurable function such that
\begin{equation}\label{eqOSKRSS10.336aFDMM}\int\limits_{\varepsilon<|z-z_0|<\varepsilon_0}K^T_{\mu}(z,z_0)\
\frac{dm(z)}{|z-z_0|^2}\ =\
o\left(\left[\log\frac{1}{\varepsilon}\right]^2\right)\qquad\forall\
z_0\in\overline D\ .\end{equation} Then Beltrami equation
(\ref{eqBeltrami}) has a multivalent solution of Dirichlet problem
(\ref{eq1002P}) for every continuous function $\varphi:E_{D}\to{\Bbb
R}$.}
\end{theorem}


\begin{remark}\label{rmOSKRSS200FDMM} Similarly, choosing in Lemma
 \ref{lemKPRS1000DMM} $\psi(t)=1/(t\log
1/t)$ instead of $\psi(t)=1/t$ we obtain that condition
(\ref{eqOSKRSS10.336aFDMM}) can be replaced by the condition
\begin{equation}\label{eqOSKRSS10.336bFDMM}
\int\limits_{\varepsilon<|z-z_0|<\varepsilon_0}\frac{K^T_{\mu}(z,z_0)\,dm(z)}{\left(|z-z_0|\
\log{\frac{1}{|z-z_0|}}\right)^2}\ =\
o\left(\left[\log\log\frac{1}{\varepsilon}\right]^2\right)\qquad\forall\
z_0\in\overline D\ .\end{equation} Here we are able to give a number
of other conditions of the logarithmic type. In particular,
condition (\ref{eqOSKRSS100dFDMM}), thanking to Theorem
\ref{thKPRS1001MM}, can be replaced by the weaker condition
\begin{equation}\label{eqOSKRSS10.336hFDMM} k^T_{z_0}(r)=O
\left(\log\frac{1}{r}\log\,\log\frac{1}{r}\right).\end{equation}
\end{remark}

Finally, by Theorem \ref{thKPRS1001MM}, applying also Theorem 3.1 in
the paper \cite{RSY}, we come to the following result.


\begin{theorem}\label{thKPRS1002MM} {\it Let $D$ be a bounded finitely connected domain in $\Bbb
C$ with nondegenerate boundary components and $\mu:{D}\to {\Bbb D}$
be a measurable function with $K_{\mu}\in L^{1}_{loc}$ such that
\begin{equation}\label{eq1009MM}\int\limits_{D\cap B(z_0,\varepsilon_0)}\Phi_{z_0}(K^T_{\mu}(z,z_0))\ dm(z)\ <\ \infty\end{equation}
for $\varepsilon_0=\varepsilon(z_0)>0$ and a nondecreasing convex
function $\Phi_{z_0}:[0,\infty)\to[0,\infty)$ with
\begin{equation}\label{eq1010MM}\int\limits_{\delta_0}^{\infty}\frac{d\tau}{\tau\Phi_{z_0}^{-1}(\tau)}\ =\ \infty\end{equation}
for $\delta_0=\delta(z_0)>\Phi_{z_0}(0)$. Then Beltrami equation
(\ref{eqBeltrami}) has a multivalent solution of Dirichlet problem
(\ref{eq1002P}) for every continuous function $\varphi:E_{D}\to{\Bbb
R}$.}
\end{theorem}


Recall that condition (\ref{eq1010MM}) is equivalent to every of
conditions (\ref{eq333Y})--(\ref{eq333D}).


\begin{corollary}\label{corKPR1003MM} {\it Let $D$ be a bounded finitely connected domain in $\Bbb
C$ with nondegenerate boundary components and $\mu:{D}\to {\Bbb D}$
be a measurable function with $K_{\mu}\in L^{1}_{loc}$ such that
\begin{equation}\label{eq1011MM} \int\limits_{D\cap B(z_0,\varepsilon_0)}e^{\alpha_0
K^T_{\mu}(z,z_0)}\ dm(z)\ <\ \infty\qquad\forall\ z_0\in\overline{D}
\end{equation}  for some
$\varepsilon_0=\varepsilon(z_0)>0$ and $\alpha_0=\alpha(z_0)>0$.

Then Beltrami equation (\ref{eqBeltrami}) has a multivalent solution
of Dirichlet problem (\ref{eq1002P}) for every continuous function
$\varphi:E_{D}\to{\Bbb R}$.}
\end{corollary}

\end{document}